\documentclass[11 pt]{article}
\usepackage{listings}
\usepackage{graphicx}
\usepackage{geometry}
\usepackage{amsmath}
\usepackage{amsfonts}
\usepackage[utf8]{inputenc}
\usepackage[english]{babel} 
\usepackage{textcomp}
\usepackage{bigints}
\usepackage{hyperref}
\usepackage{enumitem}
\usepackage{caption}
\usepackage{subcaption}
\newtheorem{theorem}{Theorem}[section]
\newtheorem{corollary}{Corollary}[theorem]
\newtheorem{remark}{Remark}[section]

\usepackage[colorinlistoftodos]{todonotes}
\RequirePackage[numbers,sort,compress]{natbib}

\begin{document}

\title{ On Robust Pseudo-Bayes Estimation for \\the Independent Non-homogeneous Set-up}
\author{
	Tuhin Majumder$^{1}$, Ayanendranath Basu$^{2\ast}$ and Abhik Ghosh$^{3\ast}$ \\\\
	$^1$ North Carolina State University, NC, USA. tmajumd@ncsu.edu\\
	$^2$ Indian Statistical Institute, Kolkata, India. ayanbasu@isical.ac.in\\
	$^3$ Indian Statistical Institute, Kolkata, India. abhik.ghosh@isical.ac.in\\
	$^\ast$Corresponding authors 
}

\maketitle

\begin{abstract}
The ordinary Bayes estimator based on the posterior density suffers from the potential problems of non-robustness 
under data contamination or outliers. In this paper, we consider the general set-up of independent but non-homogeneous (INH) observations and study a robustified pseudo-posterior based estimation for such parametric INH models. 
In particular, we focus on the $R^{(\alpha)}$-posterior developed by Ghosh and Basu (2016)\cite{Ghosh/Basu:2016}
for IID data and later extended by Ghosh and Basu (2017)\cite{Ghosh/Basu:2017} for INH set-up,
where its usefulness and desirable properties have been numerically illustrated.
In this paper, we investigate the detailed theoretical properties of this robust pseudo Bayes $R^{(\alpha)}$-posterior 
and associated $R^{(\alpha)}$-Bayes estimate under the general INH set-up with applications to fixed-design regressions. 
We derive a Bernstein-von Mises types asymptotic normality results and Laplace type asymptotic expansion
of the $R^{(\alpha)}$-posterior  as well as the asymptotic distributions of the expected $R^{(\alpha)}$-posterior estimators.
The required conditions and the asymptotic results are simplified for linear regressions with known or unknown error variance 
and logistic regression models with fixed covariates. The robustness of the $R^{(\alpha)}$-posterior and
associated estimators are theoretically examined through appropriate influence function analyses under general INH set-up;
illustrations are provided for the case of linear regression. 
A high breakdown point result is derived for the  expected $R^{(\alpha)}$-posterior estimators 
of the location parameter under a location-scale type model. 
Some interesting real life data examples illustrate possible applications.
\end{abstract}

\section{Introduction}\label{SEC:intro}

Real life data are rarely independent and identically distributed (IID); 
rather they are often non-homogeneous and can also be dependent.
Here we focus on the set-up of independent non-homogeneous (INH) observations,
which partially relaxes the IID conditions,  
with the most common example being the fixed design parametric regression set-up.
Mathematically, consider the independent random variables $X_{1}, X_{2},...,X_{n}$ which are not identically distributed; 
each $ X_{i} $ is also defined on (possibly different) measurable space $(\chi^{i}, \mathcal{B}^{i})$. 
However, we assume that there is an underlying common probability space $(\Omega,\mathcal{B}_{\Omega},P)$ such that 
$ X_{i} $ is $\mathcal{B}^{i}/\Omega$ measurable independent with respect to $P$
and its induced distribution $ G_{i}(x) $ has an absolutely continuous density $ g_{i}(x) $ with respect to 
a common dominating $\sigma$-finite measure $\lambda(dx)$, for each $i=1,2,..,n$. 
In parametric inference, we model the true unknown distributions $ G_{i}(x) $ by some appropriate
parametric families $\mathcal{F}^{i}=\{F_{i,\theta} : \theta \in \Theta \subseteq \mathbb{R}^{p}\} $, 
which are also assumed to be absolutely continuous with respect to $\lambda$ having densities 
$f_{i,\theta} $, respectively, for $i=1, \ldots, n$;
all inference are then performed using an estimate of the unknown parameter $\theta$ based on the observed data.
Although the models are different for each $i$, 
the unknown parameter $\theta$ is common in all of them leaving us with enough degrees of freedom for meaningful inference. 
In the widely used fixed-design regressions, $f_{i, \theta}$ is the density of the response variable $X_i$
for the $i$-th fixed value of the covariates, say $z_i$, and the common parameter $\theta$ consists of the regression coefficients $\beta$
along with possible variance parameters (e.g, error variance), such that $E(Y_i)=z_i^T\beta$ for $i=1, \ldots, n$.
For example, $f_{i,\theta} \equiv N(z_i^T\beta,  \sigma^2)$ 
with $\sigma^2$ being the error variance in the normal linear regression 
and $f_{i,\theta} \equiv Bernoulli(z_i^T\beta)$ density in the logistic regression model.

Among several possible approaches to do statistical inference under the above INH set-ups,
the Bayesian paradigm is one of the most popular and widely applied techniques.
One reason for this is the easy interpretibility and explainability of Bayesian procedure,
the other being its ability to incorporate the prior (to experiment) belief (often from subject knowledge)
about the final inference results. 
Considering the above INH set-up and a prior density $\pi(\theta)$ for the parameter $\theta$,
the posterior density of $\theta$ given the observed data $\underline{x}=(x_1, \ldots, x_n)$ 
is computed by the Bayes formula as
\begin{equation}
{\pi}(\theta|\underline{x})=\dfrac{\prod_{i=1}^nf_{i,\theta}(x_i)\pi(\theta)}{\int \prod_{i=1}^nf_{i,\theta}(x_i)\pi(\theta)d\theta},  
\label{EQ:post}
\end{equation}
and all the subsequent inference about $\theta$ are done based on this posterior density.
Although the resulting Bayes estimator of $\theta$ has several optimal properties, 
one major drawback of using the usual posterior density is its lack of robustness 
against data contamination or model misspecification. 
Among possible solutions to such non-robustness, a recent popular approach is to construct an appropriate pseudo-posterior
by suitably connecting this issue with the frequentist approach of robustness 
and by replacing the likelihood by some robust loss under the Bayesian paradigm;
see  \cite{Greco/etc:2008,Agostinelli/Greco:2013,Cabras/etc:2014,hv11, Ghosh/Basu:2016, Danesi/etc:2016,Atkinson/etc:2017,Nakagawa/Hashimoto:2017,Nakagawa/Hashimoto:2019, Jewson/etc:2018, Giummole/etc:2019, Futami/etc:2017, Cherief-Abdellatif/etc:2019, Gonzalez/etc:2018, Grunwald/Dawid:2004, Grunwald/VanOmmen:2017}
for several such pseudo-posteriors.
Such pseudo-posteriors are also used in other contexts like the PAC-Bayesian paradigm, 
Gibbs posterior  \cite{z99}, etc.

Unfortunately, much of the existing literature on robust Bayes inference based on pseudo-posteriors are developed only for IID data
and few special inference problems; no literature is available for our general INH set-up
except for a recent attempt by \cite{Ghosh/Basu:2017}. 
Ghosh et al.~\cite{Ghosh/Basu:2017} first studied a pseudo-posterior for general parametric models including the INH set-up
by replacing the likelihood part in (\ref{EQ:post}) with the so-called $\alpha$-likelihood (in exponential).
For our INH case, the $\alpha$-likelihood is defined as 
\begin{equation}
Q_{n}^{(\alpha)}(\theta)=\sum_{i=1}^{n}\left[\frac{1}{\alpha}f_{i,\theta}^{\alpha}(x_{i})-\frac{1}{1+\alpha}\int f_{i,\theta}^{1+ \alpha}-\frac{1}{\alpha}\right], ~~~~\alpha>0,
\label{EQ:alpha-lik}
\end{equation} 
which coincides with $\sum_{i=1}^n \log  f_{i,\theta}(x_i) - n$, the log-likelihood (up to a constant), as $\alpha\rightarrow 0$.
There are two major advantages of this particular $R^{(\alpha)}$-pseudo posterior; 
firstly, it does not require any nonparametric smoothing such as kernel density estimation unlike some other pseudo-posteriors,
and secondly, it is additive in the data so that one can update the posterior density for new observations 
without calculating it from scratch, as in the case of the usual posterior density. 
Then the corresponding pseudo-posterior, which the authors have referred to as $R^{(\alpha)}$ posterior, is defined as  
\begin{equation}\label{1}
\pi^{(\alpha)}_{R}(\theta|X_{1}, ..., X_{n})=\dfrac{\exp(Q_{n}^{(\alpha)}( \theta))\pi(\theta)}{\int \exp(Q_{n}^{(\alpha)}( \theta))\pi(\theta)d\theta}.  
\end{equation}
Note that as $\alpha\rightarrow 0$, the corresponding $R^{(\alpha)}$ posterior coincides with the ordinary Bayes posterior
and hence it provides a robust generalization at $\alpha>0$.
See \cite{Ghosh/Basu:2016,Ghosh/Basu:2017} for the motivation behind its construction,
which is based on the frequentist robust estimation procedure minimizing the density power divergence (DPD)
\cite{Basu/etc:1998,Basu/etc:2011}.
The DPD measure between two densities $f_1, f_2$ (with respect to some common dominating measure) is defined as 
\begin{eqnarray}\label{EQ:dpd}
d_\alpha(f_1,f_2) &=& \displaystyle \int  \left[f_2^{1+\alpha} - \left(1 + \frac{1}{\alpha}\right)  f_2^\alpha f_1 + 
\frac{1}{\alpha} f_1^{1+\alpha}\right]d\mu,~~~\alpha\geq 0,\\
d_0(f_1,f_2) &=& \lim\limits_{\alpha\rightarrow 0} d_\alpha(f_1,f_2) = \int f_1 \log(f_1/f_2)d\mu.
\end{eqnarray}
Note that, the  $\alpha$-likelihood in (\ref{EQ:alpha-lik}) is indeed the sum of suitably scaled negative of 
the DPD measure between $f_1=$ degenerate distribution at $x_i$ and $f_2 = f_{i,\theta}$
for $i=1, \ldots, n$.

Based on the $R^{(\alpha)}$ posterior given in (\ref{1}),
the corresponding $R^{(\alpha)}$-Bayes estimator with respect to a given loss function  $L$ is defined as 
\begin{equation*}
\widehat{\theta}_{n,\alpha}^{L}=\text{arg}\min_{t}\int L(\theta,t)\pi^{(\alpha)}_{R}(\theta|X_{1}, ..., X_{n})d\theta.
\end{equation*}
If the loss function is squared error loss, then the corresponding $R^{(\alpha)}$-Bayes estimate is 
the \textit{Expected $R^{(\alpha)}$ Posterior Estimator} (ERPE) defined as 
\begin{equation*}
\widehat{\theta}_{n,\alpha}^{E}=\int \theta\pi^{(\alpha)}_{R}(\theta|X_{1}, ..., X_{n})d\theta.
\end{equation*}
The usefulness of these robust  $R^{(\alpha)}$-Bayes estimators has already been illustrated numerically for IID data 
as well as for the INH set-up including the fixed design linear and logistic regressions.
However, theoretical properties of the $R^{(\alpha)}$-posterior has been studied in detail only for the IID cases,
which cover the exponential consistency results, concentration rate via Bernstein von-Mises type results, asymptotic normality of the corresponding ERPE, influence function analyses as well as prior robustness \cite{Ghosh/Basu:2016,Ghosh/Basu:2017}.
For INH cases, only exponential  consistency results are derived in \cite{Ghosh/Basu:2017} using appropriate assumptions 
based on the concepts of merging of distributions in probability, prior negligibility and 
existence of uniform exponential consistent tests.

In this paper, we further study the theoretical properties of the $R^{(\alpha)}$-posterior,
in both efficiency and robustness aspects, under the general INH set-up
with applications to several fixed design regression models. 
In particular, we derive the Bernstein-von Mises type asymptotic normality result  along with the asymptotic expansion for the generalized $R^{(\alpha)}$-posterior. 
The robustness properties of the corresponding $R^{(\alpha)}$-Bayes estimators are also studied 
theoretically via influence function analysis and a particular breakdown point result.
Simplifications are provided for linear and logistic regression models with applications. 

The rest of the paper is organized as follows. 
In Section \ref{SEC:BVM}, we present a Bernstein-von Mises type asymptotic normality result 
for the general INH set-up, whereas the asymptotic expansion result is provided in Section \ref{SEC:Expansion}
along with possible applications in approximating $R^{(\alpha)}$-Bayes estimators. 
Section \ref{SEC:IF} illustrates the robustness properties of the generalized $R^{(\alpha)}$ Bayes estimators 
through influence function analyses. 
In Section \ref{SEC:BP}, we present the asymptotic breakdown point for the ERPE for a location-scale type INH model. 
Additionally, in each section, the necessary conditions are simplified 
for the particular examples of the linear regression models with known or unknown error variance 
and the logistic regression. Finally, the paper ends with some concluding remarks in Section \ref{SEC:conclusions}.
For brevity in presentation, all proofs and technical calculations are moved to the Appendix.

\section{Asymptotic Normality of the Generalized \texorpdfstring{$R^{(\alpha)}$}{a}-Posterior}\label{SEC:BVM}

Consider the independent non-homogeneous set-up as in Section \ref{SEC:intro} and 
let us now study the asymptotic normality of our $R^{(\alpha)}$-posterior and the associated $R^{(\alpha)}$-Bayes estimates.
In this regard, let us recall  the minimum DPD estimator (MDPDE), say $\widehat{\theta}_{n,\alpha}$, under the INH set-up
which is obtained by minimizing $Q_{n}^{(\alpha)}(\theta)$ in (\ref{EQ:alpha-lik}) with respect to $\theta\in\Theta$;
corresponding estimating equation is then given by 
\begin{equation}
\nabla Q_{n}^{(\alpha)}(\theta)=0, ~~~\mbox{or }~~ 
\sum_{i=1}^{n}\nabla V_{i}(X_{i},\theta) = 0,
\end{equation}
where $V_{i}(x,\theta)=\int f_{i,\theta}^{1+ \alpha}-(1+\frac{1}{\alpha})f_{i,\theta}^{\alpha}(x)$
and $\nabla$ denotes the first order derivative with respect to $\theta$.
Note that, $Q_{n}^{(\alpha)}(\theta)=-\frac{1}{1+\alpha} \sum_{i=1}^{n}V_{i}(X_{i},\theta)$. 
We also assume that there exists a best fitting parameter value, say $\theta_g$, 
defined as the minimizer of the DPD measure between $g_i$ and $f_{i,\theta}$ with respect to $\theta$,
which is independent of the index $i$. It holds, for example, 
if all the true densities $g_{i}$ belong to the respective model family with $g_{i}=f_{i}(\cdot ; \theta_g)$ for all $i$s.
It is known that the MDPDE $\widehat{\theta}_{n,\alpha}$ is $\sqrt{n}$-consistent for $\theta_{g}$ 
and also asymptotically normal \cite{Ghosh/Basu:2013}, 
a property which we will exploit in the following along with its equivalence with 
our pseudo-Bayes ERPE.
Further, we will use the  notation
\begin{eqnarray}
\Psi_{n,\alpha}(\theta)&=&\dfrac{1}{1+\alpha}\dfrac{1}{n}\sum\limits_{i=1}^{n}E_{g_{i}}[\nabla^2 V_{i}(X_{i},\theta)],
\nonumber\\ 
\Omega_{n,\alpha}(\theta)&=&\dfrac{1}{n}\sum\limits_{i=1}^{n}Var_{g_{i}}[\nabla V_{i}(X_{i},\theta)]
\nonumber\\
\mbox{and }~~~\widehat{\Psi}_{n,\alpha}(\theta)&=&
\dfrac{1}{1+\alpha}\dfrac{1}{n}\sum\limits_{i=1}^{n}\nabla^2 V_{i}(X_{i},\theta)
=-\dfrac{1}{n}\nabla^2 Q_{n}^{(\alpha)}(\theta),
\nonumber
\end{eqnarray}
where $\nabla^2$ represents the second order derivative matrix with respect to $\theta$.


We start with a Bernstein-von Mises (BVM) type asymptotic normality result for the $R^{(\alpha)}$-posterior densities 
which provides its contraction rate under the general INH set-up.
For this purpose, we need to make the following assumptions:

\begin{itemize}
	\item[(E1)] The Conditions (A1)-(A6) of Ghosh and Basu \cite{Ghosh/Basu:2013} hold,
	implying the consistency of the MDPDE $\widehat{\theta}_{n,\alpha}$ for $\theta_{g}$.

\item[(E2)] For any $\delta > 0$, with probability tending to one, we have
\begin{equation}
\sup_{||\theta-\theta_{g}||>\delta}\frac{1}{n}\left[Q_{n}^{(\alpha)}(\theta)-Q_{n}^{(\alpha)}(\theta_{g})\right] < - \epsilon
\label{EQ:cond_E2}
\end{equation}
for some $\epsilon>0$ and for all sufficiently large $n$.

\item[(E3)] $0<|\Psi_{n,\alpha}(\theta)|< \infty$ for all sufficiently large $n$ and  for all $\theta \in \Theta$ such that $\lVert \theta \rVert_2 < \infty$,
where $|A|$ denotes the determinant of a matrix $A$.
\end{itemize}

The above conditions are not very hard to examine in different scenarios;
they are verified later for some common examples.
Under these assumptions, we have the following theorem.

\begin{theorem}
Suppose Assumptions (E1)--(E3) hold under the general INH set-up and 
let $\pi(\theta)$ be any prior which is positive and continuous at $\theta_{g}$. 
Then, with probability tending to one,
\begin{equation}
\lim\limits_{n \rightarrow \infty} \int \left|\pi^{*R}_{n}(t)-\left (\frac{|\Psi_{n,\alpha}(\theta_{g})|}{2 \pi}\right)^{\frac{p}{2}} e^{-\frac{1}{2}t'\Psi_{n,\alpha}(\theta_{g})t}\right| =0,
\label{EQ:BVM_genINH}
\end{equation}
where $\pi^{*R}_{n}(t)$ is the $R^{(\alpha)}$-posterior density of $t=\sqrt{n}(\theta- \widehat{\theta}_{n,\alpha})$ given the data $X_{1},X_{2},..., X_{n}$. Further, the result (\ref{EQ:BVM_genINH}) also holds with $\Psi_{n,\alpha}(\theta_{g})$ being replaced by 
$\widehat{\Psi}_{n,\alpha}(\widehat{\theta}_{n,\alpha})$. 
\label{THM:BVM_genINH}
\end{theorem}

Note that the asymptotic normality result of the $R^{(\alpha)} $-posterior under the IID set-up, as proved in \cite{Ghosh/Basu:2016},
now follows directly from the above theorem as a special case. 
However, unlike the original Bernstein-von Mises theorem, 
here we have proved convergence with probability tending to one (convergence in probability).
Indeed, if we assume the conditions under which $\widehat{\theta}_{n,\alpha}$ is strongly consistent for $\theta^{g}$ instead of Assumption (E1)
and if (E2) holds with probability one, the convergence result in Theorem \ref{THM:BVM_genINH} becomes almost sure convergence;
this can be directly observed from the proof of the theorem  given in Appendix \ref{APP:BVM_INH_pf}.

Our next theorem gives the asymptotic properties of the ERPE under our general INH set up.

\begin{theorem}
In addition to the Assumptions of Theorem \ref{THM:BVM_genINH},  assume that the prior $\pi(\theta)$ has finite expectation. 
Then the expected $R^{(\alpha)}$-posterior estimator (ERPE) $\widehat{\theta}_{n,\alpha}^{E}$ satisfies the following. 
\begin{itemize}
	\item[(a)] $\sqrt{n}\left(\widehat{\theta}_{n,\alpha}^{E} - \widehat{\theta}_{n,\alpha}\right) \xrightarrow[]{\mathcal{P}} 0$ as $n \rightarrow \infty$.
	
	\item[(b)] If, further, $\Sigma_n(\theta^{g})^{-\frac{1}{2}}\sqrt{n}\left(\widehat{\theta}_{n,\alpha}-\theta^{g}\right)
	\xrightarrow[]{\mathcal{D}} N_p(0,I_p)$ for some  positive definite matrix $\Sigma_n(\theta^{g})$, then $\Sigma_n(\theta^{g})^{-\frac{1}{2}}\sqrt{n}\left(\widehat{\theta}_{n,\alpha}^E - \theta^{g}\right)\xrightarrow[]{\mathcal{D}} N_p(0,I_p)$.
\end{itemize}
\label{THM:Asymp_ERPE}
\end{theorem}

Note that Part (a) of the above Theorem \ref{THM:Asymp_ERPE} automatically implies 
its Part (b) in view of the Slutsky theorem; Part (a) is proved later in Section \ref{SEC:Expansion}.
Additionally, it was proved in \cite{Ghosh/Basu:2013} that, under conditions (A1)--(A7) there, $\Omega_{n,\alpha}^{-\frac{1}{2}}\Psi_{n,\alpha}\left[\sqrt{n}(\widehat{\theta}_{n,\alpha}-\theta^{g})\right]\xrightarrow[]{\mathcal{D}} N_p(0,I_p)$,
which in turn then gives the asymptotic distribution of the ERPE as 
\begin{eqnarray}
\Omega_{n,\alpha}^{-\frac{1}{2}}(\theta^{g})\Psi_{n,\alpha}(\theta^{g})^{-\frac{1}{2}}\sqrt{n}
\left(\widehat{\theta}_{n,\alpha}^E - \theta^{g}\right)\xrightarrow[]{\mathcal{D}} N_p(0,I_p).
\label{EQ:Asymp_ERPE}
\end{eqnarray}

We will now simplify the above asymptotic results for some common fixed-design regression models described in the introduction. 
In this regard, let us consider the following assumptions on the fixed design points (covariates)
$z_i=(z_{i1}, \ldots, z_{ip})^T$ for $i=1, \ldots, n$. 

\begin{itemize}
\item[(R1)] For all $ j, k, l =1, \ldots, p$, we have
\begin{equation}
\displaystyle{\sup _{n \geq 1 }} \displaystyle{\max _{1 \leq i \leq n }} |z_{ij}|=O(1), 
~~~~
 \displaystyle{\max _{1 \leq i \leq n }} |z_{ij}||z_{ik}|=O(1),
~~~\mbox{and}~~~
\dfrac{1}{n}\sum_{i=1}^{n}|z_{ij}z_{ik}z_{il}|=O(1).
\end{equation}

\item[(R2)] The $n\times p$ design matrix $Z=[z_1 ~ z_2~ \cdots z_p]^T$ satisfies
$\displaystyle{\inf_{n\geq 1 }}~[\text{min eigenvalue of } n^{-1}(Z^TZ)] > 0$,
implying $Z$ has full column rank, and
$\displaystyle{\max _{1 \leq i \leq n }} [z_i^T(Z^TZ)^{-1}z_i] = O(n^{-1})$.
\item[(R3)] For any parameter value $\theta$, the matrix $\Psi_{n,\alpha}(\theta)$ satisfies
$\displaystyle{\inf_{n\geq 1 }}~[\text{min eigenvalue of } \Psi_n(\theta)] > 0$.
\end{itemize}

\bigskip\noindent
\textbf{Example \ref{SEC:BVM}.1 [Linear Regression model with known error variance]:}\\
Let us first consider the simple linear regression model $x_i=\beta z_i+\epsilon_i$ for $i=1, \ldots, n$,
where $\epsilon_i$s are IID each distributed as $N(0,\sigma^2)$ and $z_i$s are given fixed covariates. 
For the time being, we assume the error variance $\sigma^2$ to be known 
(the unknown $\sigma$ case is considered in the next example).
As mentioned in the introduction, this regression model is a special case of our INH set-up
with $\theta=\beta$ and $f_{i,\theta}$ being the conditional density of $x_i$, i.e., $N(z_i^T\beta, \sigma^2)$ density.
Using (\ref{EQ:alpha-lik}), the $\alpha$-likelihood function for this model has a simple form given by 
\begin{equation}
Q_{n}^{(\alpha)}(\beta)=\frac{1}{(2\pi)^{\frac{\alpha}{2}} \sigma^\alpha}
\sum_{i=1}^{n}\left[\frac{1}{\alpha} e^{-\frac{\alpha(x_i - z_i^T\beta)^2}{2\sigma^2}} -\frac{1}{(1+\alpha)^{\frac{3}{2}}}-\frac{1}{\alpha}\right], ~~~~\alpha>0.
\label{EQ:alpha-lik_LRM1}
\end{equation} 
Given any prior $\pi$ on the regression parameter $\beta$, one can then easily obtain the $R^{(\alpha)}$-posterior
using the $\alpha$-likelihood in (\ref{EQ:alpha-lik_LRM1}) and then compute the $R^{(\alpha)}$-Bayes estimate,
namely the ERPE of $\beta$ through an appropriate iterative numerical approach;
its performance has been empirically studied in \cite{Ghosh/Basu:2017}.
Here, we will verify Assumptions (E1)--(E3) in order to derive the theoretical (asymptotic) properties of 
the ERPE, say $\widehat{\beta}_{n,\alpha}^E$ under this linear regression model (with known $\sigma$).

Firstly, we note that (E1) has already been verified in \cite{Ghosh/Basu:2013}  when the fixed covariate values 
satisfy Assumptions (R1)--(R2). Also, it is easy to see that, here 
$\Psi_{n,\alpha}(\beta)=\zeta_{\alpha}n^{-1}{(Z^TZ)}$
with $\zeta_\alpha= (2\pi)^{-\frac{\alpha}{2}} \sigma^{-(\alpha+2)}(1+\alpha)^{-\frac{3}{2}}>0$ 
(and hence (R2) and (R3) are equivalent).
Then, it follows that $|\Psi_{n,\alpha}(\beta)|=\zeta_{\alpha}^pn^{-p}|{Z^TZ}| > 0$ 
since $Z$ is of full column rank for all $n$ by (R2). 
Since Assumption (R1) further guarantees that $|\Psi_{n,\alpha}(\beta)| < \infty $,
Condition (E3) is also implied by the simplified Assumptions (R1)--(R2). 
We will now show that the Conditions (R1)--(R2) also imply Condition (E2).
If $\beta^{g}$ denote the best fitting (true) regression parameter,
some simple algebra establishes the relation 
$$
\dfrac{1}{n}\left[Q_n(\beta)-Q_n(\beta^{g})\right]=\dfrac{1}{\alpha(\sqrt{2\pi}\sigma)^{\alpha}}\dfrac{1}{n} \sum_{i=1}^{n}\Bigg[e^{-\frac{\alpha(x_i-z_i^T\beta)^2}{2\sigma^2}}-e^{-\frac{\alpha(x_i-z_i^T\beta^{g})^2}{2\sigma^2}}\Bigg].
$$
Note that, $\sigma^{-1}{(x_i-z_i^T\beta^{g})}$ are IID ${N}(0,1)$. So, by the SLLN,  we have
$$
\dfrac{1}{n}\sum_{i=1}^{n}e^{-\dfrac{\alpha(x_i-z_i^T\beta^{g})^2}{2\sigma^2}}
\xrightarrow[]{a.s} \frac{1}{\sqrt{1+\alpha}}.
$$
For the other term, note that $\sigma^{-1}{(x_i-z_i^T\beta)}$ are independent but not identically distributed under true model
$g_i\equiv N(z_i^T\beta^g, \sigma^2)$,
but $e^{-\frac{\alpha(x_i-z_i^T\beta)^2}{2\sigma^2}}\leq 1$; hence we can apply Kolmogorov SLLN to get
$$
\begin{aligned}
&\dfrac{1}{n}\sum_{i=1}^{n}e^{-\dfrac{\alpha(x_i-z_i^T\beta)^2}{2\sigma^2}}
-\dfrac{1}{n}\sum_{i=1}^{n}E_{g_i}\Bigg[e^{-\dfrac{\alpha(X_i-z_i^T\beta)^2}{2\sigma^2}}\Bigg]\xrightarrow[]{a.s}0,
\\
\text{i.e}\quad & \dfrac{1}{n}\sum_{i=1}^{n}e^{-\dfrac{\alpha(x_i-z_i^T\beta)^2}{2\sigma^2}}- \dfrac{1}{n\sqrt{1+\alpha}}\sum_{i=1}^{n} e^{-\dfrac{\alpha(z_i^T\beta^{g}-z_i^T\beta)^2}{(1+\alpha)2\sigma^2}}\xrightarrow[]{a.s}0.
\end{aligned}
$$
Combining, we get
$$
\dfrac{1}{n}\left[Q_n(\beta)-Q_n(\beta^{g})\right]
=\dfrac{1}{\alpha(\sqrt{2\pi}\sigma)^{\alpha}}\dfrac{1}{n\sqrt{1+\alpha}}\sum_{i=1}^{n}\Bigg[ e^{-\dfrac{\alpha(z_i^T\beta^{g}-z_i^T\beta)^2}{(1+\alpha)2\sigma^2}}-1\Bigg]+U_n,
$$
where $U_n\xrightarrow[]{a.s}0.$
Now using mean value theorem on $e^{-z}$, we get
$$
e^{-\dfrac{\alpha(z_i^T\beta^{g}-z_i^T\beta)^2}{(1+\alpha)2\sigma^2}}-1=-\Bigg[\dfrac{\alpha(z_i^T\beta^{g}-z_i^T\beta)^2}{(1+\alpha)2\sigma^2}\Bigg]e^{-K_i},\quad K_i \in \Big(0,\dfrac{\alpha(z_i^T\beta^{g}-z_i^T\beta)^2}{(1+\alpha)2\sigma^2}\Big).
$$
Note that,
$
K_i \leq \dfrac{\alpha}{(1+\alpha)2\sigma^2} (\beta-\beta^{g})^Tz_iz_i^T(\beta-\beta^{g}) \leq \dfrac{\alpha}{(1+\alpha)2\sigma^2} (z_i^Tz_i)||\beta-\beta^{g}||^2.
$
But (R1) implies $z_i^Tz_i=\sum_{j=1}^{p}z^2_{ij} \leq Cp$, for some $C>0$. 
So, considering the region $\delta<\lVert \beta-\beta^{g}\rVert<2\delta$, we get
$$
K_i< \dfrac{\alpha}{(1+\alpha)2\sigma^2}Cp(2\delta)^2= \dfrac{2Cp\alpha\delta^2}{(1+\alpha)\sigma^2}=K_0 (\text{say}).
$$
Hence,\vspace{-.6cm}
$$
\dfrac{1}{n}\left[Q_n(\beta)-Q_n(\beta^{g})\right] < -\dfrac{e^{-K_0}} {2\sigma^2(\sqrt{2\pi}\sigma)^{\alpha}\sqrt{1+\alpha}}(\beta-\beta^{g})^T[n^{-1}(Z^TZ)](\beta-\beta^{g})+U_n.
$$
Now, if $\lambda_p^{(n)}$ denotes the smallest eigenvalue of $n^{-1}(Z^TZ)$, 
we know $\lambda_p^{(n)}> \zeta$ for some $\zeta>0$ under (R2). Then, for 
$\lVert \beta-\beta^{g}\rVert>\delta$, we get
$
(\beta^{g}-\beta)^T[n^{-1}(Z^TZ)](\beta^{g}-\beta)
> \zeta (\beta^{g}-\beta)^T(\beta^{g}-\beta) > \zeta \delta^2.
$ 
Using the fact that $U_n \xrightarrow[]{a.s}0$, we finally get that, 
in the region $\delta<||\beta-\beta^{g}||<2\delta$, 
$$
\dfrac{1}{n}\left[Q_n(\beta)-Q_n(\beta^{g})\right]< 
-\dfrac{e^{-K_0}} {4\sigma^2(\sqrt{2\pi}\sigma)^{\alpha}\sqrt{1+\alpha}}\zeta \delta^2,
~~~\mbox{ for all sufficiently large } n. 
$$
Finally, for any $\beta \in \mathbb{R}^p$ satisfying $\lVert  \beta - \beta^{g} \rVert \geq 2\delta$, 
we can find $\beta^{*}\in \mathbb{R}^p$ such that $\beta-\beta^{g}=\eta(\beta^{*}-\beta^{g})$ with $\eta>1$ 
and $\lVert  \beta^{*}-\beta^{g} \rVert\in (\delta, 2\delta)$. 
Then,
$(z_i^T\beta-z_i^T\beta^{g})^2=\eta^2(z_i^T\beta^*-z_i^T\beta^{g})^2\geq (z_i^T\beta^*-z_i^T\beta^{g})^2$
and hence, for all sufficiently large $n$, we get
$$
\dfrac{1}{n}(Q_n(\beta)-Q_n(\beta^{g})) \leq \dfrac{1}{n}(Q_n(\beta^*)-Q_n(\beta^{g}))<-\dfrac{e^{-K_0}} {4\sigma^2(\sqrt{2\pi}\sigma)^{\alpha}\sqrt{1+\alpha}}\zeta \delta^2.
$$
This implies that (E2) holds with 
$\epsilon=\dfrac{e^{-K_0}} {4\sigma^2(\sqrt{2\pi}\sigma)^{\alpha}\sqrt{1+\alpha}}\zeta \delta^2$.

Therefore, in this fixed-design linear regression set-up with known $\sigma^2$, 
Conditions (E1)-(E3) are implied by (R1)--(R2), and hence the  BVM type asymptotic normality result of 
Theorem \ref{THM:BVM_genINH} holds for the $R^{(\alpha)}$-posterior of $\beta$ 
under the easily verifiable  Assumptions (R1)--(R2). 
The asymptotic distribution of the corresponding ERPE $\widehat{\beta}_{n,\alpha}^E$ of $\beta$ 
can then be obtained from (\ref{EQ:Asymp_ERPE}) with the calculations of $\Omega_{n,\alpha}$;
by some simple algebra $\Omega_{n,\alpha}(\beta)=\zeta_{\alpha}n^{-1}{(Z^TZ)}$,
which leads to 
\begin{eqnarray}
(Z^TZ)^{\frac{1}{2}}\left(\widehat{\beta}_{n,\alpha}^E - \beta^{g}\right)
\xrightarrow[]{\mathcal{D}} N_p(0, \upsilon_\alpha^{(\beta)} I_p).
\label{EQ:Asymp_ERPE_LRM1}
\end{eqnarray}
where $\upsilon_\alpha^{(\beta)}= \frac{\zeta_{2\alpha}}{\zeta_\alpha^2} 
=\sigma^2 \left(1 + \frac{\alpha^2}{1+2\alpha}\right)^{\frac{3}{2}}$.
This result can now be also used to compute asymptotic relative efficiency of 
our $R^{(\alpha)}$-Bayes estimator of $\beta$, and many more properties;
further applications will be given in Section \ref{SEC:BP}.
\hfill{$\square$}

\bigskip\noindent
\textbf{Example \ref{SEC:BVM}.2 [Linear Regression model with unknown error variance]:}\\
Let us now consider the simple linear regression model as in Example \ref{SEC:BVM}.1, 
but with unknown error variance $\sigma^2$. Then, it also belongs to our INH set-up,
but now with $\theta=(\beta, \sigma)$ and $f_{i,\theta}\equiv N(z_i^T\beta, \sigma^2)$.
Interestingly, in this case also the $\alpha$-likelihood function has the same form as in (\ref{EQ:Asymp_ERPE_LRM1})
but as a function of both $\beta$ and $\sigma$, 
and hence we will denote it as $Q_{n}^{(\alpha)}(\beta, \sigma)$ instead of $Q_{n}^{(\alpha)}(\beta)$.
We show that Assumptions (R1)--(R2) are again sufficient to prove the asymptotic normality 
of the $R^{(\alpha)}$-posterior and the corresponding ERPEs of $(\beta, \sigma)$ 
even in this case with unknown $\sigma$.

Firstly, it is easy to see that Conditions (E1) and (E3) again hold under (R1)--(R2)
by arguments in line with those given in Example \ref{SEC:BVM}.1.
To verify (E2), note that, here we have
$$
\begin{aligned}
\dfrac{1}{n}\left[Q_n(\beta,\sigma)-Q_n(\beta^{g},\sigma_g)\right]&=\dfrac{1}{\alpha(\sqrt{2\pi}\sigma)^{\alpha}}\dfrac{1}{n} \sum_{i=1}^{n}e^{-\frac{\alpha(Y_i-z_i^T\beta)^2}{2\sigma^2}}-\dfrac{1}{(1+\alpha)^{3/2}(\sqrt{2\pi}\sigma)^{\alpha}}
\\& -\dfrac{1}{\alpha(\sqrt{2\pi}\sigma_g)^{\alpha}}\dfrac{1}{n} \sum_{i=1}^{n}e^{-\frac{\alpha(Y_i-z_i^T\beta^g)^2}{2\sigma_g^2}}+\dfrac{1}{(1+\alpha)^{3/2}(\sqrt{2\pi}\sigma_g)^{\alpha}}.
\end{aligned}
$$
Again, by a similar application of the Kolmogorov SLLN as in Example \ref{SEC:BVM}.1 (for known $\sigma$),
one  can observe that
\begin{equation}
\begin{aligned}
\dfrac{1}{n}\left[Q_n(\beta,\sigma)-Q_n(\beta^{g},\sigma_g)\right]&=\dfrac{1}{\alpha(\sqrt{2\pi}\sigma)^{\alpha}}\dfrac{1}{\sqrt{1+\frac{\alpha\sigma_g^2}{\sigma^2}}}\dfrac{1}{n} \sum_{i=1}^{n}\left[e^{-\frac{\alpha(z_i^T\beta-z_i^T\beta^g)^2}{(1+\alpha)2\sigma^2}}-1\right]\\
&-\dfrac{1}{\alpha(\sqrt{2\pi})^{\alpha}}\left[\dfrac{\alpha}{(1+\alpha)^{3/2}\sigma^{\alpha}}+\dfrac{1}{(1+\alpha)^{3/2}\sigma_g^{\alpha}}-\dfrac{1}{\sigma^{\alpha}\sqrt{1+\frac{\alpha\sigma_g^2}{\sigma^2}}}\right]+U_n,
\end{aligned}
\label{EQ:1}
\end{equation}
where $U_n\xrightarrow{a.s}0$.
Now, we need to prove (\ref{EQ:cond_E2}) when $\lVert (\beta^g,\sigma^g)- (\beta,\sigma)\rVert>\delta$. 
Note that, if $\sigma=\sigma_g$, the derivation boils down to the known $\sigma$ case, verified in Example \ref{SEC:BVM}.1. 
So, here we carry on our analysis for the case $\sigma \neq \sigma_g$.

Clearly, $e^{-\frac{\alpha(z_i^T\beta-z_i^T\beta^g)^2}{2\sigma^2}}-1 \leq 0$ 
and hence
the first term in the right hand side of (\ref{EQ:1}),  along with $U_n$, is negative for all sufficiently large $n$. 
To handle the second term there, we note that
$$
\dfrac{\alpha}{(1+\alpha)^{3/2}\sigma^{\alpha}}+\dfrac{1}{(1+\alpha)^{3/2}\sigma_g^{\alpha}}-\dfrac{1}{\sigma^{\alpha}\sqrt{1+\frac{\alpha\sigma_g^2}{\sigma^2}}}
=\dfrac{1}{\sigma^{\alpha}(1+\alpha)^{3/2}\sqrt{1+\alpha v}}\left[f(v) - f(1)\right],
$$
where $v=\frac{\sigma_g^2}{\sigma^2}$ and $f(v)=\alpha\sqrt{1+\alpha v}+v^{-\alpha/2}\sqrt{1+\alpha v}$. 
Now, by simple calculations, one can verify that $f'(v)>0$ for $v>1$ and $f'(v)<0$ for $v<1$, 
which implies that  $f$ is strictly decreasing for $v<1$ and strictly increasing for $v>1$. 
Hence, $f$ has the unique minima at $v=1$, i.e $f(v)>f(1)$ for $v \neq 1$, i.e for $\sigma \neq \sigma_g$. 
Hence, if $|\sigma-\sigma_g|>\delta$, then either $v>\delta_1>1$ or $v<\delta_2<1$
for some $\delta_1, \delta_2>0$. Then, the second term in the right hand side of (\ref{EQ:1})
becomes strictly less than $-C$ for some $C>0$, implying Condition (E2).

Therefore, under the present fixed-design linear regression set-up with unknown $\sigma^2$ also, 
Assumptions (R1)--(R2) imply the  asymptotic normality of the $R^{(\alpha)}$-posterior of $(\beta, \sigma)$,
as in Theorem \ref{THM:BVM_genINH}, and that of their respective ERPEs 
$\widehat{\beta}_{n,\alpha}^E$ and $\widehat{\sigma}_{n,\alpha}^E$, as in (\ref{EQ:Asymp_ERPE}).
The calculations of the matrices $\Psi_{n,\alpha}$ and $\Omega_{n,\alpha}$ can be done in a straightforward manner 
(also available in \cite{Ghosh/Basu:2013}) which shows that 
the asymptotic distribution of $\widehat{\beta}_{n,\alpha}^E$ is the same as in (\ref{EQ:Asymp_ERPE_LRM1})
and that of $\widehat{\sigma}_{n,\alpha}^E$ is independent of $\widehat{\beta}_{n,\alpha}^E$ with
\begin{eqnarray}
\sqrt{n}\left(\widehat{\sigma}_{n,\alpha}^E - \sigma_{g}\right)
\xrightarrow[]{\mathcal{D}} N_p(0, \upsilon_\alpha^{(\sigma)} I_p).
\label{EQ:Asymp_ERPE_LRM2}
\end{eqnarray}
where $\upsilon_\alpha^e 
= \frac{\sigma^2}{(2+\alpha^2)^2}
\left[2(1+2\alpha^2)\left(1+ \frac{\alpha^2}{1+2\alpha}\right)^{\frac{5}{2}} - \alpha^2(1+\alpha)^2\right]$.
The asymptotic relative efficiency of both the ERPEs are independent of the values of true parameter and the design matrix;
the values for some particular $\alpha$ are provided in Table \ref{TAB:ARE_LRM}. Clearly, 
there is no significant loss in asymptotic efficiency at small values of $\alpha>0$.
\hfill{$\square$}

\begin{table}[!h]
	\centering
	\caption{The asymptotic relative efficiencies of the ERPEs of $(\beta, \sigma)$ under the fixed-design linear regression model with unknown error variance}
	\begin{tabular}{c c c c c c c c c c c} 
\hline
		$\alpha$ & 0 & 0.01 & 0.02 & 0.05 & 0.10 & 0.15 & 0.25 & 0.50 & 0.75 &	1.00 \\ \hline
		\hline
		$\beta$ & 100 & 99.99 & 99.94 & 99.66 & 98.76 & 97.46 & 94.06 & 83.81 & 73.76 & 64.95 \\
		$\sigma$ & 100 & 99.97 & 99.88 & 99.32 & 97.56 & 95.05 & 88.84 & 73.06 & 61.53 & 54.11 \\
		\hline
	\end{tabular}
\label{TAB:ARE_LRM}
\end{table}

\bigskip\noindent
\textbf{Example \ref{SEC:BVM}.3 [Fixed-design Logistic Regression model]:}\\
Our next example is the highly popular logistic regression model for binary response data.
Given fixed covariate values  $z_i$s, here, the binary response $x_i$ is assumed to follow a Bernoulli distribution
with success probability $\pi_{i}(\beta) =\dfrac{e^{z_i^{T}\beta}}{1+e^{z_i^{T}\beta}}$. 
Clearly, it is also a special case of our INH set-up
with $\theta=\beta$ and 
$$
f_{i,\theta}(x_i)=f_{i,\beta}(x_i)=\dfrac{e^{x_i (z_i^{T}\beta)}}{1+e^{z_i^{T}\beta}}, 
~~~\mbox{for }~x_i=0, 1.
$$ 
Using (\ref{EQ:alpha-lik}), the corresponding $\alpha$-likelihood function for 
this fixed-design logistic regression model can be simplified as 
\begin{equation}
Q_{n}^{(\alpha)}(\beta)=\sum_{i=1}^{n}\left[\dfrac{e^{\alpha x_i (z_i^{T}\beta)}}{\alpha (1+e^{z_i^{T}\beta})^\alpha} 
-\frac{1+e^{(1+\alpha) (z_i^{T}\beta)}}{(1+\alpha)(1+e^{z_i^{T}\beta})^{1+\alpha}}-\frac{1}{\alpha}\right], ~~~~\alpha>0.
\label{EQ:alpha-lik_Log}
\end{equation} 
Now, for any appropriate prior $\pi$ on the regression coefficient $\beta$, 
we can obtain the $R^{(\alpha)}$-posterior from (\ref{EQ:alpha-lik_LRM1}) using (\ref{EQ:alpha-lik_Log}).
The $R^{(\alpha)}$-Bayes estimates of $\beta$ under any loss function can then be computed numerically;
the particular case of ERPE is discussed in \cite{Ghosh/Basu:2017}.

Let us now verify study their theoretical (asymptotic) properties  with verification of the required  Assumptions (E1)--(E3).
Note that, in this context, some straightforward algebra leads to
$$
\Psi_{n,\alpha}(\beta)=\dfrac{1}{n}\sum_{i=1}^{n} \dfrac{e^{z_i^{T}\beta}(e^{\alpha z_i^{T}\beta}+e^{z_i^{T}\beta})}{(1+e^{z_i^{T}\beta})^{3+\alpha}}z_iz_i^{T},
~~~
\Omega_{n,\alpha}(\beta)=\dfrac{1}{n}\sum_{i=1}^{n} \dfrac{e^{z_i^{T}\beta}(e^{\alpha z_i^{T}\beta}+e^{z_i^{T}\beta})^2}{(1+e^{z_i^{T}\beta})^{4+2\alpha}}z_iz_i^{T}.
$$
Here, Assumption (R3) directly implies (E3) with the above form of $\Psi_{n,\alpha}$.
Indeed, Assumptions (R1) and (R3) are also sufficient to imply Conditions (E1) and (E2);
the proofs are a bit involved and hence are moved to Appendix \ref{APP:logistic_verification} for brevity. 
Therefore, the easily verifiable Assumptions (R1) and (R3) imply the BVM type asymptotic normality result,
as in Theorem \ref{THM:BVM_genINH}, for the $R^{(\alpha)}$-posterior of $\beta$ 
under the fixed-design logistic regression models. 
Accordingly, using (\ref{EQ:Asymp_ERPE}) with the value of $\Omega_{n,\alpha}$ as given above,
the asymptotic distribution of the ERPE of $\beta$ can be obtained in a straightforward manner. 
\hfill{$\square$}

\section{Asymptotic Expansions and Approximate Computation of the $R^{(\alpha)}$-Bayes estimates }
\label{SEC:Expansion}

As in the usual Bayes theory, the computation of the $R^{(\alpha)}$-Bayes estimates
also needs the evaluation of the integrals of the form 
$\int h(\theta) \exp(Q_n^{(\alpha)}(\theta)) \pi(\theta) d\theta$ 
for some given function $h(\theta)$ and prior $\pi(\theta)$.
Since these integrals often do not have an expression in explicit closed form
and numerical integration is mostly infeasible beyond dimension one or two, 
one prefers a suitable large sample approximation for these integrals .
In usual Bayes inference, it is often achieved via Laplace expansion method,
a powerful technique for good large sample approximations. 

Let us now study such approximations for our ERPE under the general INH set-up. 
Here, we denote $B_{\delta}(\theta)$ to be an open ball of radius $\delta>0$ around $\theta\in\Theta$,
and $\nabla_{j_1\ldots j_d}$ to be the $d$-th order derivative with respect to $\theta_{j_1}, \ldots, \theta_{j_d}$ 
for any $1 \leq j_1,...,j_d \leq p$ and $d\geq 1$.
Also recall that the MDPDE obtained by minimizing the $\alpha$-likelihood $Q_n^{(\alpha)}(\theta)$
is denoted as $\widehat{\theta}_{n,\alpha}$; see Section \ref{SEC:BVM}.
Then, we have the following theorem for a large sample expansion of the $R^{(\alpha)}$-posterior based integrals
under INH set-up; the proof follows directly from  \cite{Kass/etc:1990} and is hence omitted for brevity.

\begin{theorem}
Consider the general INH set-up and assume that the corresponding $\alpha$-likelihood $Q_n^{(\alpha)}(\theta)$, 
given in (\ref{EQ:alpha-lik}), 	is six times continuously differentiable function on $\Theta$ for each $n\geq1$.
Also consider a four times continuously differentiable function $q(\theta)$ on $\Theta\subseteq \mathbb{R}^p$.
Suppose there exists positive number $\epsilon$, $M$ and $\eta$ and an integer $n_0$ such that 
the followings are satisfied for all $n \geq n_0$.
 \begin{enumerate}\label{analytical}
     \item[(B1)] With $0 \leq d \leq 6$, we have $n^{-1}|\nabla_{j_1...j_d}Q_n^{(\alpha)}(\theta)| < M$ 
     for all $1 \leq j_1,...,j_d \leq p$ and $\theta \in B_{\delta}(\theta)$.
     \item[(B2)] $\left|-n^{-1}\nabla^2 Q_n^{(\alpha)}(\widehat{\theta}_{n,\alpha})\right|> \eta$.
     \item[(B3)] For all $\delta\in (0, \epsilon)$, $B_{\delta}(\widehat{\theta}_{n,\alpha}) \subset \Theta$ and 
$$
\left|-\nabla^2 Q_n^{(\alpha)}(\widehat{\theta}_{n,\alpha})\right|^{1/2}
\int_{\Theta\setminus B_{\delta}(\widehat{\theta}_{n,\alpha})} q(\theta)
\exp\left[ Q_n^{(\alpha)}({\theta}) - Q_n^{(\alpha)}(\widehat{\theta}_{n,\alpha})\right]d\theta
=O\left(n^{-2}\right).
$$
\end{enumerate}
Then, we have the following asymptotic expansion 
\begin{eqnarray}
\int q(\theta)\exp[Q_n^{(\alpha)}(\theta)]d\theta =
q(\widehat{\theta}_{n,\alpha})\exp[{Q_n^{(\alpha)}(\widehat{\theta}_{n,\alpha})}](2\pi)^{\frac{p}{2}}
\left|-\nabla^2 Q_n^{(\alpha)}(\widehat{\theta}_{n,\alpha})\right|^{-1/2} \Big[1+O(n^{-1})\Big].
\label{EQ:Expansion_INH}
\end{eqnarray}
\label{THM:Expansion_INH}
\end{theorem}

\begin{corollary}
Under the assumptions of Theorem \ref{THM:Expansion_INH} with $q(\theta)=h(\theta)\pi(\theta)$
and $q(\theta)=\pi(\theta)$, 
we have the following expansion of the expectation of $h(\theta)$ under the $R^{(\alpha)}$-posterior distribution of $\theta$
as given by
$$
\frac{\int h(\theta)\exp[Q_n^{(\alpha)}(\theta)]\pi(\theta)d\theta}{\int \exp[Q_n^{(\alpha)}(\theta)]\pi(\theta)d\theta}
=h(\widehat{\theta}_{n,\alpha})\Big[1+O(n^{-1})\Big].
$$
In particular, the ERPE satisfies  $\widehat{\theta}_{n,\alpha}^{E}= \widehat{\theta}_{n,\alpha}\Big[1+O(n^{-1})\Big]$,
implying Result (a) of Theorem \ref{THM:Asymp_ERPE}. 
\end{corollary}

\bigskip
Note that, verifying the Assumptions (B1)--(B3) for any particular examples of INH set-up,
we can use the expansion (\ref{EQ:Expansion_INH}) to approximate any $R^{(\alpha)}$-Bayes estimator
including the ERPE; we will present some examples below. 
However, in some cases, it might be easier to verify a stronger version (B3)$^\ast$ (given below) of Assumption (B3),
which is similar in  spirit to Assumption (E2); see also \cite{Kass/etc:1990}.
\begin{enumerate}
	\item[(B3)$^\ast$] For all $\delta\in (0, \epsilon)$, $B_{\delta}(\widehat{\theta}_{n,\alpha}) \subset \Theta$ and 
	$\displaystyle{\limsup_{n \rightarrow \infty }} 
	\displaystyle{\sup_{\theta\in \Theta\setminus B_{\delta}(\widehat{\theta}_{n,\alpha}) }}
	\frac{1}{n}\left[Q_{n}^{(\alpha)}(\theta)-Q_{n}^{(\alpha)}(\widehat{\theta}_{n,\alpha})\right]<0.$
\end{enumerate}

\bigskip\noindent
\textbf{Example \ref{SEC:Expansion}.1 [Fixed-design Linear Regression model]:}\\
Let us again consider the fixed-design linear regression model as described in Examples \ref{SEC:BVM}.1, \ref{SEC:BVM}.2.
We will now verify the Assumptions (B1)--(B3) so that the Laplace expansion in (\ref{EQ:Expansion_INH})
can be used for approximate computation of our $R^{(\alpha)}$-Bayes estimators of parameter in this regression model. 

Note that, whether the error variance $\sigma^2$ is known or unknown,
Assumption (B1) is straightforward to verify and (B2) is directly implied by Condition (R2).
Now, we argue that Conditions (R1)--(R2) also imply Assumption (B3)$^\ast$ and hence (B3). 
For this purpose, let us first assume $\sigma$ to be known (so that $\theta=\beta$)
and the true density belongs to the model family (so that $g_i=f_{i, \beta^g}$ for all $i$).
Then, after some algebra, we get
\begin{eqnarray}
&&\frac{1}{n}\left[Q_{n}^{(\alpha)}(\theta)-Q_{n}^{(\alpha)}(\widehat{\beta}_{n,\alpha})\right]
\nonumber\\
&=& \dfrac{1}{\sqrt{1+\alpha}(\sqrt{2\pi} \sigma)^{\alpha}}\dfrac{1}{n}\sum_{i=1}^{n}
\left[e^{-\frac{\alpha(\beta-\beta^g)^Tz_iz_i^T(\beta-\beta^g)}{2\sigma^2(\alpha+1)}}
- e^{-\frac{\alpha(\widehat{\beta}_{n,\alpha}-\beta^g)^Tz_iz_i^T(\widehat{\beta}_{n,\alpha}-\beta^g)}{2\sigma^2(\alpha+1)}}
\right] \nonumber\\ 
&=&\dfrac{1}{\alpha\sqrt{1+\alpha}(\sqrt{2\pi} \sigma)^{\alpha}}\dfrac{1}{n}\sum_{i=1}^{n}
\left[e^{-\frac{\alpha(\beta-\beta^g)^Tz_iz_i^T(\beta-\beta^g)}{2\sigma^2(\alpha+1)}}-1\right]+ o_p(1)
\end{eqnarray}
since $\widehat{\beta}_{n,\alpha} \rightarrow \beta^g$ under Assumptions (R1)--(R2) \cite{Ghosh/Basu:2013}. 
Now, we can follow the calculations in the line of Example \ref{SEC:BVM}.1 to get  some $K>0$, which satisfies
$$
\frac{1}{n}\left[Q_{n}^{(\alpha)}(\theta)-Q_{n}^{(\alpha)}(\widehat{\beta}_{n,\alpha})\right] \leq 
-\dfrac{ e^{-K}}{2\sigma^2(1+\alpha)^{1/2}(\sqrt{2\pi} \sigma)^{\alpha}}(\beta-\beta^g)^T[n^{-1}(Z^TZ)](\beta-\beta^g)+o_p(1).
$$
If $\lambda_n^{(p)}$ denote the smallest eigenvalue of the matrix $[n^{-1}(Z^TZ)]$, 
for all sufficiently large $n$ on the region $\Theta\setminus B_{\delta}(\widehat{\beta}_{n,\alpha})$, we have
$$
(\beta-\beta^g)^T[n^{-1}(Z^TZ)](\beta-\beta^g)
\geq \lambda_n^{(p)}(\beta-\beta^g)^T(\beta-\beta^g)
>\lambda_n^{(p)} (\lVert \beta-\widehat{\beta}_{n,\alpha} \rVert- \lVert \widehat{\beta}_{n,\alpha} -\beta^g \rVert)^2 
> \eta_0 \dfrac{\delta^2}{4},
$$
due to the consistency of $\widehat{\beta}_{n,\alpha}$. 
Therefore, the required inequality of (B3)$^\ast$ holds and hence Assumption (B3).

For the case of unknown error variance also, we can similarly modify the calculations of Example \ref{SEC:BVM}.2
to show that (B3)$^\ast$ again holds under the same Assumptions (R1)--(R2).
\hfill{$\square$}

\bigskip\noindent
\textbf{Example \ref{SEC:Expansion}.2 [Fixed-design Logistic Regression model]:}\\
Here, we want to verify assumptions (B1)--(B3) for logistic regression model described in Example \ref{SEC:BVM}.3. (B1) is straightforward to show, and (B2) follows from the condition (R3). Only challenge remains is to verify (B3)$^\ast$. But this verification is simiilar to the verification of (E1)--(E3), as done in appendix. Using the fact that $\hat{\beta}_{n,\alpha}\rightarrow \beta^g$, we can say that 
$$
\begin{aligned}
\frac{1}{n}\left[Q_{n}^{(\alpha)}(\beta)-Q_{n}^{(\alpha)}(\widehat{\beta}_{n,\alpha})\right] &= \frac{1}{n}\left[Q_{n}^{(\alpha)}(\beta)-Q_{n}^{(\alpha)}(\beta^g)\right]+
\frac{1}{n}\left[Q_{n}^{(\alpha)}(\beta^g)-Q_{n}^{(\alpha)}(\widehat{\beta}_{n,\alpha})\right]\\
&\leq 
-\epsilon+o_p(1)\leq -\dfrac{\epsilon}{2}
\end{aligned}
$$
for sufficiently large $n$ and for $\Theta\setminus B_{\delta}(\widehat{\beta}_{n,\alpha})$.
\section{Local Robustness under Data Contamination}\label{SEC:IF}

The influence function (IF) is a classical frequentist measure of (local) robustness of an estimator against data contamination,
which measures the asymptotic bias caused by an infinitesimal amount contamination at a distant outlying point. 
It has been used in the Bayesian context by \cite{Ghosh/Basu:2016} to study the theoretical robustness properties 
of the $R^{(\alpha)}$-Bayes estimators as well as the whole $R^{(\alpha)}$-posteriors density.
Although IF is initially proposed and mostly discussed in the context of IID data,
it has also been extended for the general INH set-up by \cite{Huber:1983,Ghosh/Basu:2013} and others. 
Here, we now use the extended definitions of the IF to theoretically justify the robustness
of the $R^{(\alpha)}$-posterior based inference under the general INH set-up.

Consider the INH set-up and associated notations as in Section \ref{SEC:intro}.
Let us denote the collections of the true and model distribution functions as
$\mathcal{G}=\{G_1, \ldots, G_n \}$ and  $\mathcal{F}_{n,\theta}=\{F_{1,\theta},...,F_{n,\theta}\}$, respectively. 
In order to define the IF for the INH set-up, one considers the contaminated distributions $G_{i,\epsilon}=(1-\epsilon)G_{i}+\epsilon \wedge_{t_{i}}$, 
where $\wedge_{t_{i}}$ is the degenerate distribution at the point of contamination $t_{i}$ for $i=1, \ldots, n$. 
Then, the IF of any statistical  functional $T(\mathcal{G})$, associated with the estimator under consideration,
can be defined with contamination either  at a fixed $i_{0}$-th direction or in all the directions, respectively, as 
\begin{eqnarray}
IF_{i_{0}}(t_{i_{0}}; T,\mathcal{G}) &=&
\frac{\partial}{\partial \epsilon}T_{\alpha}(G_{1},...,G_{i_{0},\epsilon},...,G_{n})|_{\epsilon=0},
\nonumber\\
IF(t_{1},..,t_{n};T,\mathcal{G}) &=& 
\frac{\partial}{\partial\epsilon}T_{\alpha}(G_{1,\epsilon},...,G_{n,\epsilon})|_{\epsilon=0}.
\nonumber
\end{eqnarray}
If these IFs are bounded in the contaminations points, the corresponding estimator is expected to be robust 
whereas unbounded IF indicates the non-robust nature of the underlying estimator.

Since the IF is defined in terms of statistical functionals, in our present context,
we first need to define a functional corresponding to the $\alpha$-likelihood $Q_{n}^{(\alpha)}(\theta)$ as
\begin{eqnarray}
Q^{(\alpha)}(\theta; \mathcal{G},\mathcal{F}_{n,\theta})
=\sum\limits_{i=1}^{n}\Bigg[ \dfrac{1}{\alpha}\displaystyle\int  f^{\alpha}_{i,\theta}(x)dG_{i}(x)
-\dfrac{1}{1+\alpha}\displaystyle\int f_{i,\theta}^{1+\alpha}(x)dx - \frac{1}{\alpha}\Bigg]
=\sum\limits_{i=1}^{n}Q^{(\alpha)}_{i}(\theta; G_{i},F_{i,\theta}).
\label{EQ:alpka-lik_func}
\end{eqnarray}
Following the ideas of \cite{Ghosh/Basu:2016}, we can then write the generalized $R^{(\alpha)}$-posterior density 
as a functional of the true data generating distributions $\mathcal{G}$ and the unknown parameter $\theta$ as follows:
\begin{equation}
\pi_{\alpha}(\theta, \mathcal{G})=\dfrac{\pi(\theta)e^{Q^{(\alpha)}(\theta; \mathcal{G},\mathcal{F}_{n,\theta})}}{\int \pi(\theta)e^{Q^{(\alpha)}(\theta; \mathcal{G},\mathcal{F}_{n,\theta})}d\theta}.
\label{EQ:R-alpha-post_func}
\end{equation} 
For fixed sample size $n$, given a loss function $L(\cdot,\cdot)$,
the corresponding $R^{\alpha}$-Bayes functional can then be defined as
\begin{equation}
T_{n}^{(\alpha)L}(\mathcal{G})=\text{arg}\min_{t}\dfrac{\int L(\theta,t)\pi(\theta)e^{Q^{(\alpha)}(\theta; \mathcal{G},\mathcal{F}_{n,\theta})}d\theta}{\int \pi(\theta)e^{Q^{(\alpha)}(\theta; \mathcal{G},\mathcal{F}_{n,\theta})}d\theta}.
\label{EQ:R-Bayes_Est_Func}
\end{equation}
In particular, if $L$ is the squared error Loss, 
the statistical functional corresponding to the ERPE is given by
\begin{equation}
T_{n}^{(\alpha)e}(\mathcal{G})=\dfrac{\int \theta\pi(\theta)e^{Q^{(\alpha)}(\theta; \mathcal{G},\mathcal{F}_{n,\theta})}d\theta}{\int \pi(\theta)e^{Q^{(\alpha)}(\theta; \mathcal{G},\mathcal{F}_{n,\theta})}d\theta}.
\label{EQ:ERPE_Func}
\end{equation}
We now derive the IF for these ERPE functional as well as the general $R^{\alpha}$-Bayes functionals
in the following subsections.

\subsection{Influence Function of expected \texorpdfstring{$R^{(\alpha)}$}{a} Posterior}
\label{SEC:IF_ERPE}

One can easily compute the IF of the ERPE by straight-forward differentiation at a fixed sample size $n$ 
which, under contamination in only $i_{0}$-th direction, turns out to have the form
\begin{equation}
IF_{i_{0}}(t_{i_{0}},T_{n}^{(\alpha)e},\mathcal{G})=\text{Cov}_{P_{R_{n}}}(\theta,k_{i_{0},\alpha}(\theta,t_{i_{0}},g_{i_{0}})).
\end{equation} 
where $\text{Cov}_{P_{R_{n}}}$ denotes the covariance under 
the (true) $R^{(\alpha)}$-posterior functional density in (\ref{EQ:R-alpha-post_func}) and
\begin{eqnarray}
k_{i,\alpha}(\theta,t_{i},g_{i})&=&
\dfrac{1}{\alpha}\bigg[f^{\alpha}_{i,\theta}(t_{i})-\displaystyle\int f^{\alpha}_{i,\theta}g_{i}\bigg],
~~~~~\alpha>0,
\nonumber\\
k_{i,0}(\theta,t_{i},g_{i})&=& \log f_{i,\theta}(t_{i}) - \displaystyle\int g_{i} \log f_{i,\theta}.
\end{eqnarray}

Similarly, the influence function of $T_{n}^{(\alpha)e}$ with contamination at all the data points 
can be seen to have the form
\begin{equation}\label{16}
IF(t_{1},...,t_{n},T_{n}^{(\alpha)e},\mathcal{G})
=\sum\limits_{i=1}^{n}\text{Cov}_{P_{R_{n}}}(\theta,k_{i,\alpha}(\theta,t_{i},g_{i})).
\end{equation}

Although it is not straightforward to verify the boundedness of these IFs, 
for most common densities in exponential form, it is clear that the associated quantity 
$k_{i,\alpha}(\theta,t_{i},g_{i})$ is unbounded in $t_i$ at $\alpha=0$ but becomes bounded for any $\alpha>0$.
The exact values of the IF can be computed numerically to verify the same for the IF
indicating the claimed robustness of our proposed ERPE at $\alpha>0$; 
an example in the context of linear regression is given in Subsection \ref{SEC:IF_LRM} below.

\subsection{Influence Function of General \texorpdfstring{$R^{(\alpha)}$}{a} Bayes Estimator}

We now derive the influence function of the $R^{(\alpha)}$-Bayes estimator 
with respect to any general loss function $L(.,.)$, provided $L$ is twice differentiable in its second argument.
For this purpose, we note that the associated functional in (\ref{EQ:R-Bayes_Est_Func})
can alternatively be defines as the solution of the  estimating equation
\begin{center}
$\dfrac{\partial}{\partial t}\Bigg[\dfrac{\int L(\theta,t)\pi(\theta)e^{Q^{(\alpha)}(\theta; \mathcal{G},\mathcal{F}_{n,\theta})}d\theta}{\int \pi(\theta)e^{Q^{(\alpha)}(\theta; \mathcal{G},\mathcal{F}_{n,\theta})}d\theta}\Bigg]=0$.
\end{center}
Through standard differentiation, this estimating equation can be simplified as 
\begin{equation}
\displaystyle\int L'(\theta,T_{n}^{(\alpha)L}(\mathcal{G})\exp[Q^{(\alpha)}(\theta; \mathcal{G},\mathcal{F}_{n,\theta})]\pi(\theta)d\theta=0,
\label{EQ:R-Byes_Est_EE}
\end{equation}
where we denote
$L'(\theta,T_{n}^{(\alpha)L}(\mathcal{G}))$ denotes the partial derivative of $L(\theta,t)$ with respect to $t$ evaluated at  
$t=T_{n}^{(\alpha)L}(\mathcal{G})$ (Similarly we will also denote the corresponding second order derivative by $L''$).

Now, in order to derive the IF of $T_{n}^{(\alpha)L}(\mathcal{G})$ under contamination in the $i_{0}^{th}$ direction, 
we replace $G_{i_{0}}$ in the estimating equation (\ref{EQ:R-Byes_Est_EE}) by $G_{i_{0},\epsilon}$ 
and differentiate it with respect to $\epsilon$ at $\epsilon=0$. Collecting terms with some little algebra,
the resulting IF turns out to have the form 
\begin{equation}
IF_{i_{0}}(t_{i_{0}},T_{n}^{(\alpha)L},\mathcal{G})
=-\dfrac{E_{P_{R_{n}}}[L'(\theta,T_{n}^{(\alpha)L}(\mathcal{G}))k_{i_{0},\alpha}(\theta,t_{i_{0}},g_{i_{0}})]}{
	E_{P_{R_{n}}}[L''(\theta,T_{n}^{(\alpha)L}(\mathcal{G}))]}.
\nonumber
\end{equation} 

Similarly, we can replace each $G_{i}$ in the estimating equation (\ref{EQ:R-Byes_Est_EE}) by $G_{i,\epsilon}$ 
and proceed as before to get the IF of $T_{n}^{(\alpha)L}(\mathcal{G})$ under contaminate in all the distributions,
which has the form
\begin{equation} 
 IF(t_{1},...,t_{n},T_{n}^{(\alpha)L},\mathcal{G}) = -\sum\limits_{i=1}^{n}\dfrac{E_{P_{R_{n}}}[L'(\theta,T_{n}^{(\alpha)L}(\mathcal{G}))k_{i,\alpha}(\theta,t_{i},g_{i})]}{E_{P_{R_{n}}}[L''(\theta,T_{n}^{(\alpha)L}(\mathcal{G}))]}
=\sum\limits_{i=1}^{n}IF_{i}(t_{i},T_{n}^{(\alpha)L},\mathcal{G}).
\nonumber
\end{equation}

\subsection{Influence Function of overall \texorpdfstring{$R^{(\alpha)}$}{a}-Posterior Density}
\label{SEC:IF_posterior}

In Bayesian paradigm the whole posterior distribution is more important than a particular point estimator
and hence its robustness is worth investigating to study the robustness of any Bayes procedure based on it.
Ghosh and Basu \cite{Ghosh/Basu:2016} first developed some pseudo influence function (PIF) measure to investigate 
the robustness properties of a density function, having same interpretation as the usual IF,
and applied it to prove the non-robustness of the usual Bayes posterior having unbounded PIF
and also the robustness of $R^{(\alpha)}$-posterior density under the IID set-up.
Here, we now extend their definition of PIF in the context of general INH data 
to justify the claimed robustness of the corresponding generalized $R^{(\alpha)}$-posterior density functional 
defined in (\ref{EQ:R-alpha-post_func}).

Let us first consider the contamination to be at a particular $i_0$-th direction (distribution)
and measure the local change in the $R^{(\alpha)}$ Posterior density due to infinitesimal contamination as
$$
\displaystyle{\lim_{\epsilon \downarrow 0}}\dfrac{\pi_{\alpha}(\theta;G_{1},...,G_{i_{0},\epsilon},...,G_{n})-\pi_{\alpha}(\theta;\mathcal{G})}{\epsilon}
=\dfrac{\partial}{\partial \epsilon}\pi_{\alpha}(\theta;G_{1},...,G_{i_{0},\epsilon},...,G_{n}) \Big|_{\epsilon=0}.
$$
After some calculation, one can simplify this measure to the form 
 \begin{eqnarray}
\dfrac{\partial}{\partial \epsilon}\pi_{\alpha}(\theta;G_{1},...,G_{i_{0},\epsilon},...,G_{n}) \Big|_{\epsilon=0}
& =&\pi_{\alpha}(\theta;\mathcal{G})[k_{i_{0},\alpha}(\theta,t_{i_{0}},g_{i_{0}})-E_{P_{R_{n}}} \{k_{i_{0},\alpha}(\theta,t_{i_{0}},g_{i_{0}})\}]
\nonumber\\
&=&\pi_{\alpha}(\theta;\mathcal{G})\mathcal{I}_{i_{0},\alpha}(\theta , t_{i_0}, G_{i_0}),
\end{eqnarray}
where the expectation in the above expression is taken under the (true) $R^{(\alpha)}$-posterior density functional  $\pi_{\alpha}(\theta;\mathcal{G})$ in (\ref{EQ:R-alpha-post_func}) and 
$$
\mathcal{I}_{i,\alpha}(\theta , t_{i}, G_{i})=k_{i,\alpha}(\theta,t_{i},g_{i})
-E_{P_{R_{n}}} \left[k_{i,\alpha}(\theta,t_{i},g_{i})\right]
$$
provides an interpretation same as the IF. As discussed in \cite{Ghosh/Basu:2016},
we refer to it as the pseudo influence function under the INH set-up.
At a finite sample size $n$, this provides us a measure of local robustness of the whole pseudo posterior density.
In particular, if this quantity $\mathcal{I}_{i,\alpha}(\theta , t_{i}, G_{i})$ is unbounded, 
then the relative change in the $R^{(\alpha)}$-posterior density may be infinite 
under distant contamination in the $i$-th density, 
which means our inference based on this pseudo-posterior will be unstable. 
On the other hand, if $\mathcal{I}_{i,\alpha}(\theta , t_{i}, G_{i})$ is bounded, 
our inference would be more robust. 
Note that expectation of this quantity is zero under the $R^{(\alpha)}$ posterior density functional. 

Similarly, when we assume the contamination is in all distributions,
 the corresponding measure of local change in the $R^{(\alpha)}$-posterior density can be given by
\begin{eqnarray}
\displaystyle{\lim_{\epsilon \downarrow 0}}\dfrac{\pi_{\alpha}(\theta;G_{1,\epsilon},...,G_{n,\epsilon})-\pi_{\alpha}(\theta;\mathcal{G})}{\epsilon}
&=&\dfrac{\partial}{\partial \epsilon}\pi_{\alpha}(\theta;G_{1,\epsilon},...,G_{n,\epsilon}) \Big|_{\epsilon=0}
\nonumber\\
&=&
\pi_{\alpha}(\theta;\mathcal{G})\sum_{i=1}^{n}\mathcal{I}_{i,\alpha}(\theta , t_{i}, G_{i}).
\end{eqnarray}
Here the relevant quantity to measure the robustness is the total pseudo influence measure PIF
$\mathcal{I}_{\alpha}(\theta; t_{1},...,t_{n}; \mathcal{G})= \sum_{i=1}^{n}\mathcal{I}_{i,\alpha}(\theta , t_{i}, G_{i})$,
which can be verified for boundedness to study the robustness of the underlying pseudo-posterior density under INH set-up;
see Example in the next subsection.
Based on this measure $\mathcal{I}_{\alpha}(\theta; t_{1},...,t_{n}; \mathcal{G})$, 
we can also define local and global measure of sensitivity of the $R^{(\alpha)}$-posterior density 
with respect to the contamination in our INH data, respectively, 
as $\gamma_{\alpha}(t_1,...,t_n)=\sup_{\theta} \mathcal{I}_{\alpha}(\theta; t_{1},...,t_{n}; \mathcal{G})$
at any set of contamination points $\{t_1,...,t_n\}$ and $\gamma_{\alpha}^{*}=\sup\limits_{t_1,...,t_n}\gamma_{\alpha}(t_1,...,t_n)$.

An alternative interpretation of the PIF can also be provided under the general set-up
as done in \cite{Ghosh/Basu:2016} for IID cases. Let us consider the $\phi$-divergence measure 
between densities $v_{1},v_{2}$ as defined by $\rho(v_{1},v_{2})= \int \phi(\dfrac{v_{1}}{v_{2}})v_{2}$, 
where $\phi$ is a smooth convex function with bounded first and second derivative and $\phi(1)=0$. 
Then, depending on the contamination scenarios, and by usual differentiations, 
one can show that
\begin{equation}
\displaystyle{\lim_{\epsilon \downarrow 0}} \dfrac{\rho(\pi_{\alpha}(\theta;G_{1},...,G_{i,\epsilon},...,G_{n}) , \pi_{\alpha}(\theta;\mathcal{G}))}{\epsilon}=\phi'(1) E_{P_{R_{n}}}\Big[ \mathcal{I}_{i,\alpha}(\theta , t_{i}, G_{i}) \Big] = 0
\end{equation}
and
\begin{equation}
\displaystyle{\lim_{\epsilon \downarrow 0}} \dfrac{\rho(\pi_{\alpha}(\theta;G_{1,\epsilon},...,G_{n,\epsilon}) , \pi_{\alpha}(\theta;\mathcal{G}))}{\epsilon}=\phi'(1) E_{P_{R_{n}}}\Big[ \mathcal{I}_{\alpha}(\theta; t_{1},...,t_{n}; \mathcal{G})] = 0.
\end{equation}
However, if we magnifying the divergence by $\epsilon^{2}$ before taking limits, 
for contamination in all directions, we indeed get 
\begin{eqnarray}
\displaystyle{\lim_{\epsilon \downarrow 0}} \dfrac{\rho(\pi_{\alpha}(\theta;G_{1,\epsilon},...,G_{n,\epsilon}) , \pi_{\alpha}(\theta;\mathcal{G}))}{\epsilon^{2}}
&=&\dfrac{\partial^{2}}{\partial^{2} \epsilon} \rho(\pi_{\alpha}(\theta;G_{1,\epsilon},...,G_{n,\epsilon}) , \pi_{\alpha}(\theta;\mathcal{G})) \Big |_{\epsilon=0} 
\nonumber \\
&=&\phi''(1) E_{P_{R_{n}}}\Big[ \mathcal{I}_{\alpha}(\theta; t_{1},...,t_{n}; \mathcal{G})\Big]^{2} 
\nonumber \\
&=& \phi''(1) Var_{P_{R_{n}}}\Big[ \mathcal{I}_{\alpha}(\theta; t_{1},...,t_{n}; \mathcal{G})\Big].
\end{eqnarray}
This limiting value also gives a possible measure of local sensitivity 
of the $R^{(\alpha)}$-posterior density under the INH set-up, 
which we may denote by $s_{\alpha}(t_1,...,t_n)$;
A global measure of the corresponding (variance) sensitivity can be defined as 
$s_{\alpha}^*=\displaystyle{\sup_{t_1,...,t_n}} s_{\alpha}(t_1,...,t_n)$.
Similar results also hold for contamination in any one fixed-direction and is hence omitted for brevity

\subsection{Example: Linear Regression with Known Error Variance}
\label{SEC:IF_LRM}

For illustration of the pseudo IF analyses, we recall the example of fixed-design linear regression model 
with known error variance, as described in Example \ref{SEC:BVM}.1.
For simplicity, let us consider a one dimensional parameter (one covariate only); 
the cases of multiple covariates can be dealt in a similar manner.

\begin{figure}[h]
	\begin{subfigure}{.5\textwidth}
		\includegraphics[width=8cm,height=7cm]{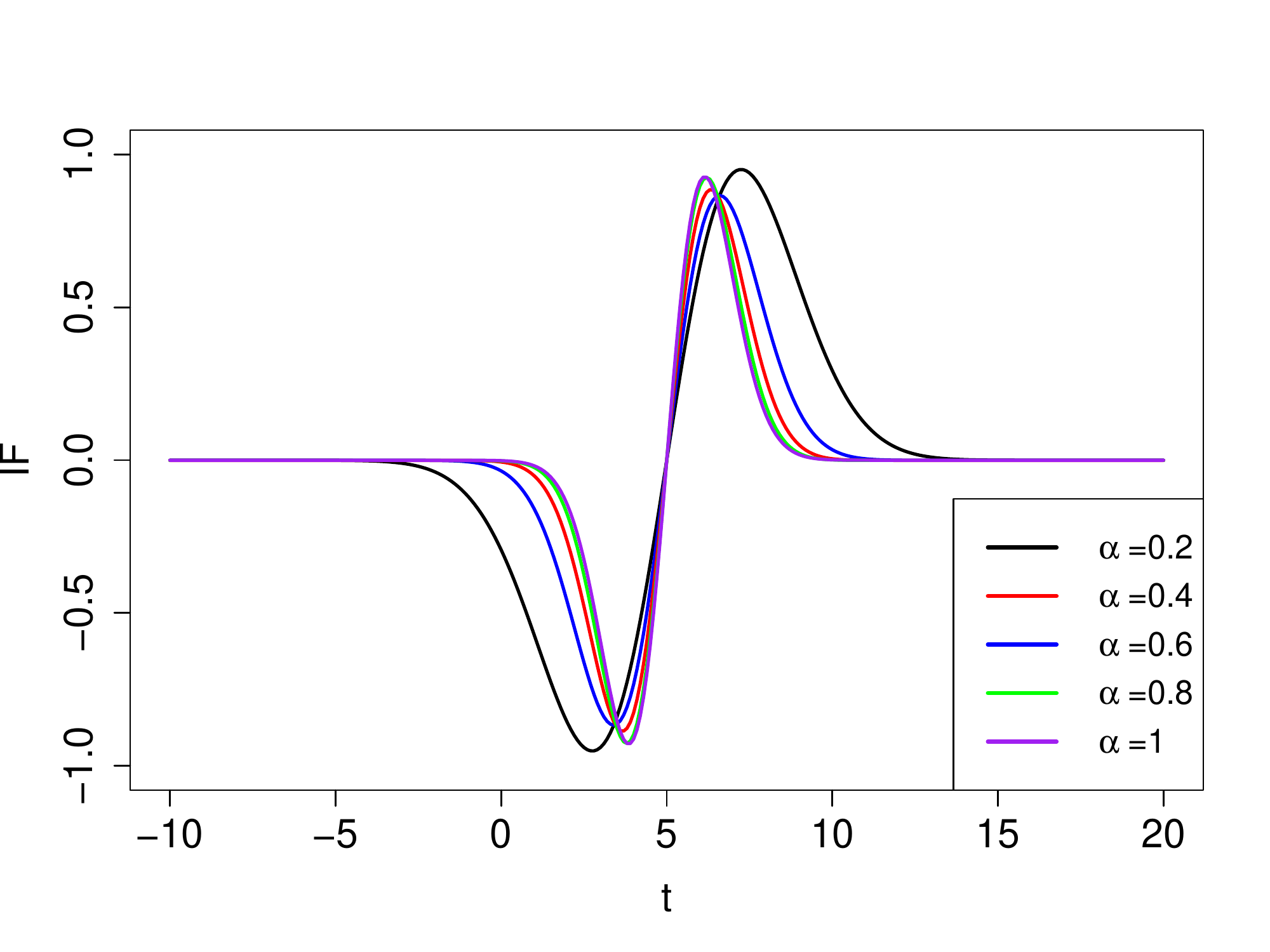}
		\caption{ n=20}
		\label{fig:sub1_1}
	\end{subfigure}
	\begin{subfigure}{.5\textwidth}
		\includegraphics[width=8 cm, height=7 cm]{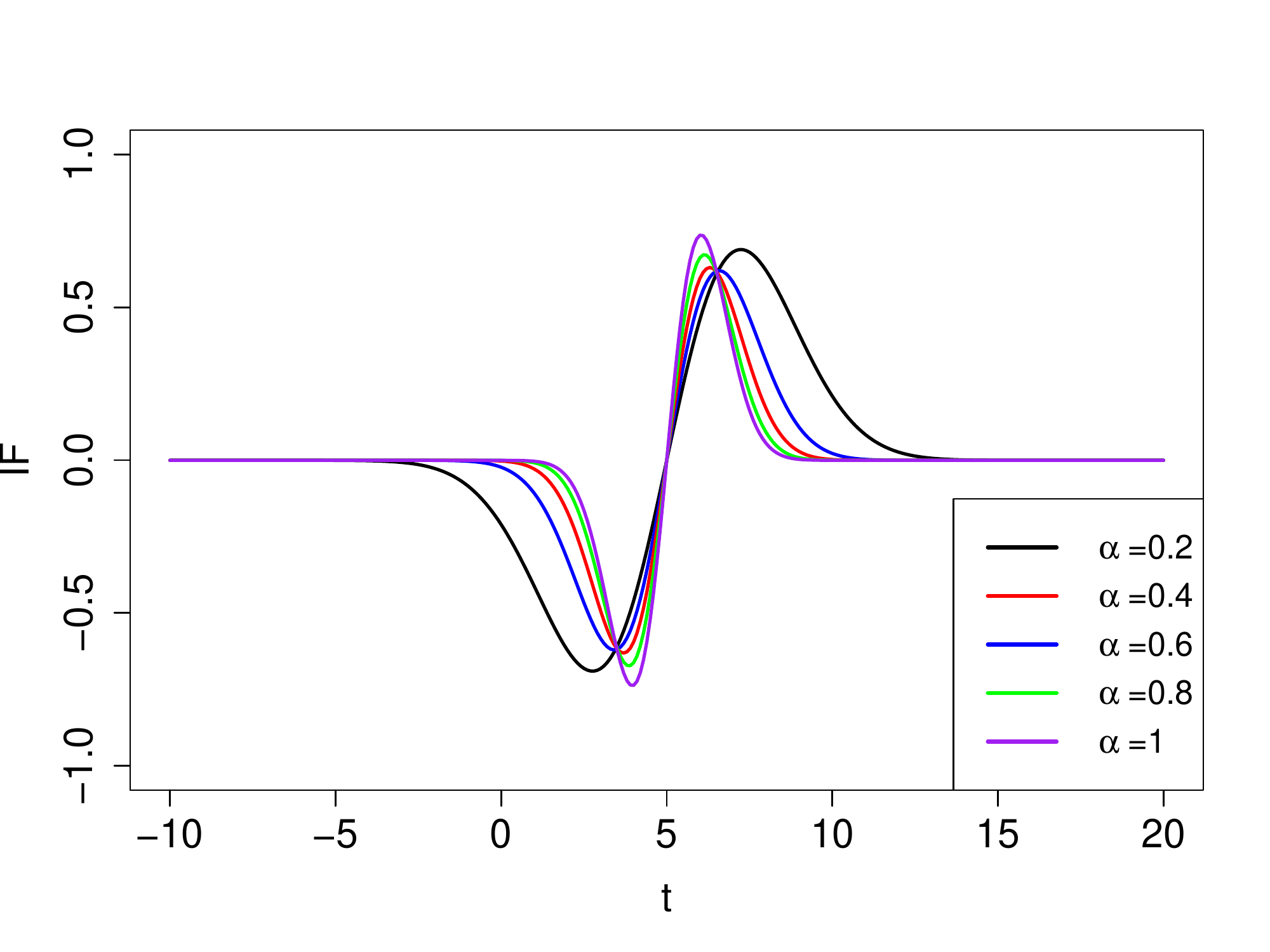}
		\caption{ n=50}
		\label{fig:sub1_2}
	\end{subfigure}
	\caption{\small Plots of the fixed sample IF of the ERPE for different $\alpha$ and sample sizes $n$}
	\label{FIG:IF_ERPE}
\end{figure}

We first study the robustness of the ERPE of the regression coefficient $\theta=\beta$ 
via the corresponding IF formula derived in Subsection \ref{SEC:IF_ERPE}.
Note that, it is not easy to get a closed form expression for the IF at $\alpha>0$, although we can compute it for $\alpha=0$ as
$$
IF(t,t,\ldots,t,T_{n}^{(0)e},\mathcal{G})=\frac{n(t-\overline{Z}\beta^g)}{n+1};
$$
here we have considered the contamination in all distributions at the points $t_{1}=t_{2}=\cdots=t_{n}=t$.
Clearly, this IF of the usual Bayes estimator under squared error loss (posterior expectation) is unbounded in $t$ 
for every fixed sample size $n$ indicating its non-robust nature even under infinitesimal contamination.
For $\alpha>0$, we numerically compute the IF of the corresponding ERPEs 
assuming contamination in all distributions at the points $t_{1}=t_{2}=\cdots=t_{n}=t$;
we consider the example with $\beta^g=5$, $z_{i}$ are generated from $N(1,1)$ and the prior to be $N(5,1)$,
and the computations are done using an appropriate importance sampling scheme. 
The resulting values of the IF are plotted over $t$ in Figure \ref{FIG:IF_ERPE}
for different values of $\alpha>0$ and two fixed sample sizes $n=20, 50$.
The boundedness of these IFs is clear from the figures indicating the claimed robustness of our ERPE at $\alpha>0$.

\begin{figure}[h]
	\begin{subfigure}{.5\textwidth}
		\includegraphics[width=8cm]{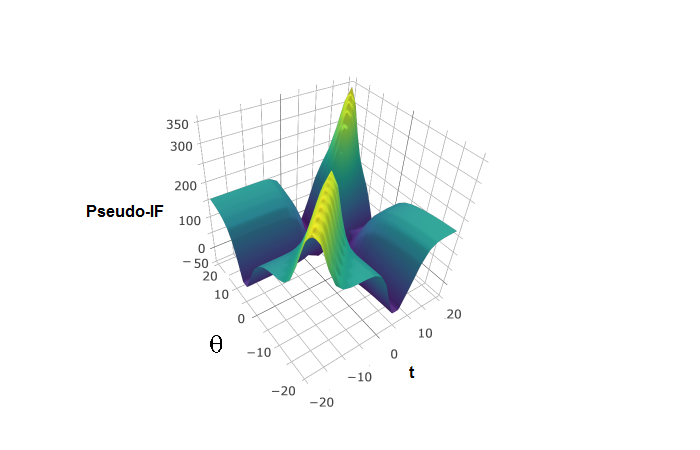}
		\caption{ $n=20$,  $\alpha= 0.1$}
		\label{fig:sub2_1}
	\end{subfigure}
	\begin{subfigure}{.5\textwidth}
		\includegraphics[width=8cm]{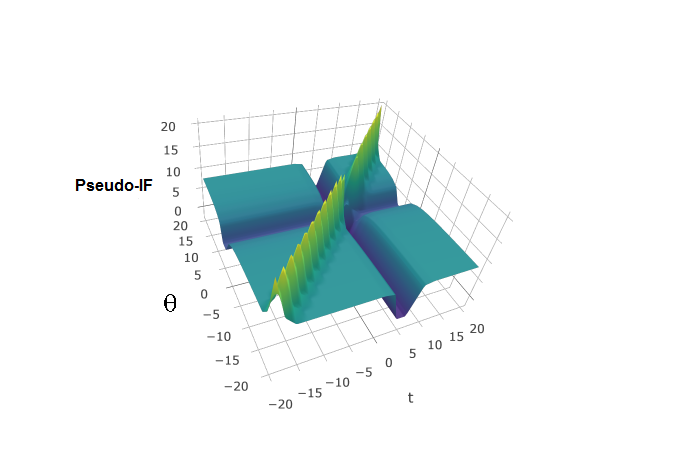}
		\caption{ $n=20$,  $\alpha= 0.8$}
		\label{fig:sub2_2}
	\end{subfigure}
	\begin{subfigure}{.5\textwidth}
		\includegraphics[width=8cm]{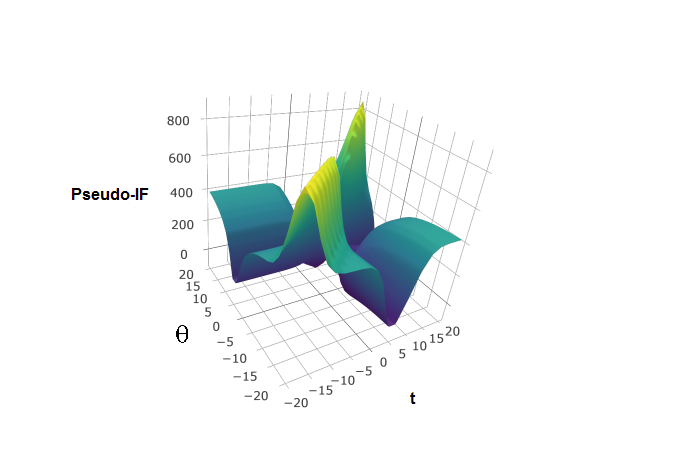}
		\caption{ $n=50$,  $\alpha= 0.1$}
		\label{fig:sub2_3}
	\end{subfigure}
	\begin{subfigure}{.5\textwidth}
		\includegraphics[width=8cm]{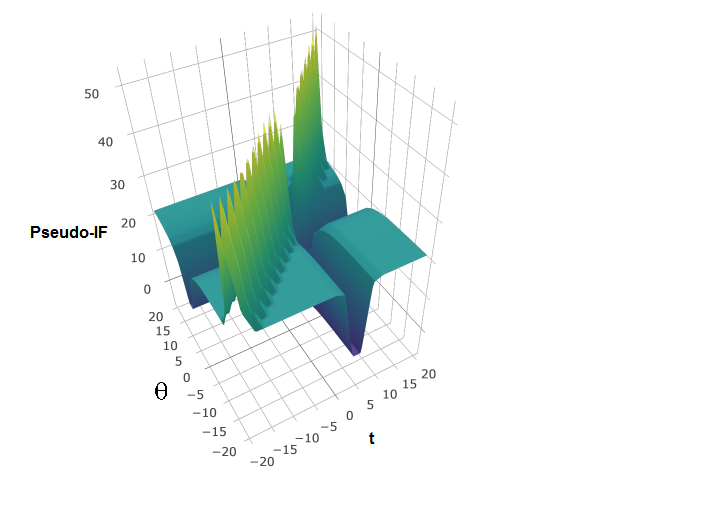}
		\caption{ $n=50$, $\alpha= 0.8$}
		\label{fig:sub2_4}
	\end{subfigure}
	\caption{\small Plots of fixed sample pseudo-influence function of $R^{(\alpha)}$ posterior density for several $\alpha$ and for different sample sizes $n$}
	\label{FIG:IF_post}
\end{figure}

Next let us study the robustness of the whole posterior density via the pseudo IF derived in Section \ref{SEC:IF_posterior}.
We again consider contamination in all distributions at the points $t_1=...=t_n=t$ and the set-ups as before
to numerically compute these pseudo IFs at $\alpha>0$, which are presented in Figure \ref{FIG:IF_post}.
At $\alpha=0$, we again can derive its explicit form as given by 
$$
\mathcal{I}_{\alpha}(\beta; t,...,t; \mathcal{G}) = \left(\beta-\beta^g\right)\sum_{i=1}^{n}\left(tz_i-\beta^gz_i^2\right),
$$
which is clearly unbounded in the contamination point $t$ indicating the non-robust nature of the usual posterior density.
However, as seen from Figure \ref{FIG:IF_post}, 
the pseudo IF of the general $R^{(\alpha)}$-posterior density is bounded at $\alpha>0$
indicating the claimed robustness of any inference based on this proposed pseudo posterior.
Further, the supremum of these pseudo IFs (sensitivities)  can be seen to decrease as $\alpha$ increases from 0.1 to 0.8
which indicates the increase in the extent of robustness with increasing values of $\alpha>0$.

\section{A Global Robustness Result for the ERPE: Breakdown Point}
\label{SEC:BP}

The global robustness of an estimator functional is measured by the (asymptotic) breakdown point. 
This frequentist property of estimators are widely used under IID data though different (but related) formulations;
see for example \cite{Hampel/etc:1986}. However, in the context of INH data, it is recently defined in \cite{Ghosh/Basu:2013}
by extending the corresponding IID idea from \cite{Simpson:1987}. 
Considering a statistical functional $T(\mathcal{G})$ at the true densities $\mathcal{G}=\{G_1, \ldots, G_n\}$ 
under the our INH set-up with a given sample size $n$, 
the breakdown point of $T$ is defined as the minimum $\epsilon>0$ such that 
there exists a sequence of contamination distributions  $\{K_{i,m}\}$ satisfying 
$$
|T((1-\epsilon)G_1 + \epsilon K_{1,m}, \ldots, (1-\epsilon)G_n + \epsilon K_{n,m}) - T(\mathcal{G})| \rightarrow \infty,
~~~\mbox{ as } m\rightarrow \infty. 
$$
Ghosh and Basu \cite{Ghosh/Basu:2013} also used this definition to prove the breakdown point of the minimum DPD functional
defined as the minimizer of the $\alpha$-likelihood functional $Q^{(\alpha)}(\theta; \mathcal{G},\mathcal{F}_{n,\theta})$ 
in (\ref{EQ:alpka-lik_func}) with respect to $\theta$ under a location-scale type INH model. 
In particular, they have considered the models families to be 
\begin{equation}
 \mathcal{F}_{i,\theta}=\Bigg\{\frac{1}{\sigma}f\Big(\frac{y-l_{i}(\mu)}{\sigma}\Big) : \theta=(\mu, \sigma) \in \Theta\Bigg\}
\label{EQ:loc-scal_family}
\end{equation} 
where $\mu$ is a location parameter and $\sigma$ is the scale parameter,
and proved that the minimum DPD functional $T_{n,\alpha}^{\mu}(\mathcal{G})$ of $\mu$ 
has the maximum asymptotic breakdown point of $0.5$ at the model with fixed scale parameter,
under certain assumptions [(BP1)--(BP3) in their paper \cite{Ghosh/Basu:2013}].

We will now apply the above idea of breakdown point in the Bayesian context to study the global robustness 
of our $R^{(\alpha)}$-Bayes estimators, in particular for the ERPE. 
The idea is to prove the asymptotic equivalence of the ERPE functional $T_{n}^{(\alpha)e}(\mathcal{G})$,
defined in (\ref{EQ:ERPE_Func}), for the location parameter $\theta=\mu$ 
under the model family in (\ref{EQ:loc-scal_family}) with known variance ($\sigma$)
with the corresponding minimum DPD functional   $T_{n,\alpha}^{\mu}(\mathcal{G})$; 
the rigorous result is presented in the following theorem for general INH models
which follows in a similar fashion as in the proofs of Section \ref{SEC:Expansion}.

\begin{theorem}
Consider the general INH set-up and assume that the conditions of Theorem \ref{THM:Expansion_INH}
hold with $\alpha$-likelihood $Q_n^{(\alpha)}({\theta}) $ replaced by the corresponding functional 
$Q^{(\alpha)}(\theta; \mathcal{G},\mathcal{F}_{n,\theta})$ and $q(\theta)$ replaced by the corresponding prior $\pi(\theta)$.
Then, an asymptotic expansion similar to (\ref{EQ:Expansion_INH}) holds for the $R^{(\alpha)}$-posterior functional 
and the ERPE functional $T_{n}^{(\alpha)e}(\mathcal{G})$, defined in (\ref{EQ:ERPE_Func}),
is asymptotically equivalent to the corresponding minimum DPD functional $T_{n,\alpha}(\mathcal{G})$,
defined as a minimizer of $Q^{(\alpha)}(\theta; \mathcal{G},\mathcal{F}_{n,\theta})$ with respect to $\theta\in \Theta$,
in the sense
$$
T_{n}^{(\alpha)e}(\mathcal{G})= T_{n,\alpha}(\mathcal{G})\Big[1+O(n^{-1})\Big].
$$
\label{THM:Expansion_Func}
\end{theorem}

Note that Theorem \ref{THM:Expansion_Func} holds for general INH models and 
not only for the particular model family given in (\ref{EQ:loc-scal_family}).
Then, the asymptotic breakdown point result of our ERPE 
under this location-scale type model (\ref{EQ:loc-scal_family})
follows directly from the corresponding results for the MDPDE given in \cite{Ghosh/Basu:2013}, 
which is presented in the subsequent theorem;
the proof is given in Appendix \ref{APP:BP} for brevity of presentation.

\begin{theorem}	
Consider the INH set-up with the model family being given by (\ref{EQ:loc-scal_family}) with fixed $\sigma$ 
(so that $\theta=\mu$) and assume that the conditions of Theorem \ref{THM:Expansion_INH}
hold with $\alpha$-likelihood $Q_n^{(\alpha)}({\theta}) $ replaced by the corresponding functional 
$Q^{(\alpha)}(\theta; \mathcal{G},\mathcal{F}_{n,\theta})$ and 
$q(\theta)$ replaced by the underlying prior $\pi(\mu)$ on $\theta=\mu$.
Further, if Assumptions (BP1)--(BP3) given in Ghosh and Basu \cite{Ghosh/Basu:2013} hold 
for the model families and the contamination densities, then the asymptotic breakdown point $\epsilon^*$ of 
the ERPE functional $\mu_{n}^{(\alpha)e}=T_{n}^{(\alpha)e}(\mathcal{G})$ of the location parameter $\mu$ at any $\alpha>0$ 
is at least $\frac{1}{2} $ at the model.
\label{THM:BP_ERPE}
\end{theorem}

Note that, this particular location scale type family (\ref{EQ:loc-scal_family})
covers the important case of the fixed-design linear regression with known error variance 
as discussed in Example \ref{SEC:BVM}.1. Hence, all the breakdown point results derived in this section 
applies to the estimators of the corresponding regression coefficient $\beta$.

\begin{remark}
Since IID set-ups are particular cases of the general INH set-up with $f_{i,\theta}=f_\theta$ for all $i$,
all the asymptotic and breakdown point results of this section also hold for the IID cases;
such results were not available in the literature of IID cases as well.
\end{remark}

\section{Conclusions}
\label{SEC:conclusions}
This paper illustrates a Bernstein-von Mises type result which provides the contraction rate of 
the robust pseudo posterior in INH set up under certain assumptions. 
It also deals with the asymptotic expansion of the robustified posterior in general. 
Not only that, we have simplified our assumptions in the well-known scenarios, 
such as linear regression and logistic regression.

The second major contribution of this paper is to develop the local robustness properties of 
the $R^{(\alpha)}$ Bayes estimators, mainly by deriving influence function and pseudo-influence functions in general. 
We have also dealt with the analysis of a global robustness measure,
namely breakdown point of the ERP functionals of the location parameter under a location-scale type INH model, 
keeping the scale parameter fixed.

In a nutshell, we think that our current work is very extensive, 
and broadens the scope of finding new research problems in this area. 
As a future work, we can think of developing similar analysis in non-independent set ups, 
such as in time series or Markov models. 
Another direction could be to develop a new Bayes hypothesis testing using the robust posterior, 
both in IID and INH setup, and generalizing all the results in existing literature for usual Bayesian inference. 
Some more applications of the proposed methodologies in generalized linear models will also be worthwhile.

\appendix
\section{Proof of Theorem \ref{THM:BVM_genINH}}\label{APP:BVM_INH_pf}

Using the form of $\pi^{(\alpha)}_{R}(\theta|X_{1},X_{2},..., X_{n})$, we get,
\begin{equation}
\begin{aligned}
&\pi^{*R}_{n}(t)=\frac{\exp(Q_{n}^{(\alpha)}(\hat{\theta}_{n,\alpha}+\frac{t}{\sqrt{n}}))\pi(\hat{\theta}_{n,\alpha}+\frac{t}{\sqrt{n}})}{\int \exp (Q_{n}^{(\alpha)}(\hat{\theta}_{n,\alpha}+\frac{t}{\sqrt{n}}))\pi(\hat{\theta}_{n,\alpha}+\frac{t}{\sqrt{n}}) dt}\\
&=c^{-1}_{n} \exp\left[Q_{n}^{(\alpha)}(\hat{\theta}_{n,\alpha}+\frac{t}{\sqrt{n}})-Q_{n}^{(\alpha)}(\hat{\theta}_{n,\alpha})\right]\pi(\hat{\theta}_{n,\alpha}+\frac{t}{\sqrt{n}}),
\end{aligned}
\end{equation} 
where $c_n$ is the required normalizing constant. Define $$
g_{n}(t)=\pi(\hat{\theta}_{n,\alpha}+\frac{t}{\sqrt{n}})\exp[Q_{n}^{(\alpha)}(\hat{\theta}_{n,\alpha}+\frac{t}{\sqrt{n}})-Q_{n}^{(\alpha)}(\hat{\theta}_{n,\alpha})]-\pi(\theta_{g})e^{-\frac{1}{2}t'\Psi_{n,\alpha}(\theta_{g})t}.
$$ 
Then to prove the first part, it is enough to show that with probability tending to one,
\begin{equation}
\int |g_{n}(t)| dt \rightarrow 0  \quad \textrm{as  } n\rightarrow \infty.
\end{equation}
\hspace*{0.25 in}Let us consider $S_{1}= \{t: ||t||>\delta_{0}\sqrt{n}\}$ and $S_{2}= \{t: ||t|| \leq \delta_{0}\sqrt{n}\}$. We shall show that $\int_{S_{i}}|g_{n}(t)| dt \rightarrow 0$, as $n\rightarrow \infty$ for $i=1,2$. 
From now, for notational simplicity, we omit the superscript $\alpha$ in $Q_{n}^{(\alpha)}$ and $\hat{\theta}_{n,\alpha}$;  and just denote it by $Q_{n}$ and $\hat{\theta}_{n}$ respectively. So, by definition,
\begin{center}
$\nabla Q_{n}(\hat{\theta}_{n}) =0 .$
\end{center} 
\hspace*{0.25 in}Now, as per our notation, $-\dfrac{1}{n}\bigtriangledown^{2}Q_{n}(\theta)=\hat{\Psi}_{n}(\theta)$.
By assumption (A6) of \cite{Ghosh/Basu:2013}, and applying a generalized version of Khinchin's weak law, it follows that $\forall j,k$, $(\hat{\Psi}_{n}(\theta_{g}))_{jk}-(\Psi_{n,\alpha}(\theta_{g}))_{jk} \rightarrow 0$ with probability tending to one. Now,
$$
\begin{aligned}
Q_{n}(\hat{\theta}_{n}+\frac{t}{\sqrt{n}})-Q_{n}(\hat{\theta}_{n})&= \dfrac{1}{2n}t'[\nabla^{2} Q_{n}(\hat{\theta}_{n})]t+\dfrac{1}{6n\sqrt{n}}\sum_{j,k,l}t_{j}t_{k}t_{l}\bigtriangledown_{jkl}Q_{n}(\theta_{n}')\\
& =-\dfrac{1}{2}t'\hat{\Psi}_{n}(\hat{\theta}_{n})t+ R_{n}(t).
\end{aligned}
$$
Note that 
$$
\begin{aligned}
|R_{n}(t)| & \leq C\dfrac{1}{\sqrt{n}}\sum_{j,k,l}t_{j}t_{k}t_{l}\left|\dfrac{1}{n}\sum_{i} \bigtriangledown_{jkl}V_{i}(X_{i},\theta_{n}')\right| \\ &\leq C\dfrac{1}{\sqrt{n}}\sum_{j,k,l}t_{j}t_{k}t_{l}\left(\dfrac{1}{n}\sum_{i}\left| \bigtriangledown_{jkl}V_{i}(X_{i},\theta_{n}')\right|\right)\\
&\leq C\dfrac{1}{\sqrt{n}}\sum_{j,k,l}t_{j}t_{k}t_{l}\left(\dfrac{1}{n}\sum_{i}M^{(i)}_{jkl}(X_{i}) \right).
\end{aligned}
$$ 

By assumption (A5) , $|\dfrac{1}{n}\sum_{i}M^{(i)}_{jkl}(X_{i}) | < 2m_{jkl} < \infty $ with prob tending to 1 for sufficiently large $n$. Thus for a fixed $t$, $R_{n}(t) \rightarrow 0$ as $n \rightarrow \infty$. Now for $1 \leq j,k \leq p$, $$
\begin{aligned}
(\hat{\Psi}_{n}(\hat{\theta}_{n}))_{jk} &=-\dfrac{1}{n}\bigtriangledown_{jk} Q_{n}(\hat{\theta}_{n})\\& =-\dfrac{1}{n}\bigtriangledown_{jk} Q_{n}(\theta_{g})+\dfrac{1}{1+ \alpha}\sum_{l} (\hat{\theta}_{n}^{(l)}-\theta^{(l)}_{g})\dfrac{1}{n}\sum_{i=1}^{n} \bigtriangledown_{jkl}V_{i}(X_{i},\theta_{n}'')\\
&=(\hat{\Psi}_{n}(\theta_{g}))_{jk}+R_{1n}(t).
\end{aligned}
$$

Now, by similar method as we used to show $R_{n}(t) \rightarrow 0$,  
$$
|R_{1n}(t)| \leq C_{1}\sum_{l}| (\hat{\theta}_{n}^{(l)}-\theta^{(l)}_{g})|\left(\dfrac{1}{n}\sum_{i}M^{(i)}_{jkl}(X_{i}) \right).
$$
Since $|\dfrac{1}{n}\sum_{i}M^{(i)}_{jkl}(X_{i}) | < \infty $ with prob tending to 1 for sufficiently large $n$ and $ \hat{\theta}_{n} $ is consistent for $\theta_{g}$, we conclude that $R_{1n}(t)=o_{p}(1)$.
So $$
\begin{aligned}
\hat{\Psi}_{n}(\hat{\theta}_{n})=\Psi_{n,\alpha}(\theta_{g})+(\hat{\Psi}_{n}(\theta_{g})-\Psi_{n,\alpha}(\theta_{g})) + o_{p}(1)=\Psi_{n,\alpha}(\theta_{g})+o_{p}(1),
\end{aligned}
$$
since we have shown $(\hat{\Psi}_{n}(\theta))_{jk}-(\Psi_{n,\alpha}(\theta))_{jk} \rightarrow 0
$ with probability tending to one.
Hence for a fixed $t$, 
$$
\begin{aligned}
|g_{n}(t) |&=|\pi(\hat{\theta}_{n}+\frac{t}{\sqrt{n}})\exp[Q_{n}(\hat{\theta}_{n}+\frac{t}{\sqrt{n}})-Q_{n}(\hat{\theta}_{n})]-\pi(\theta_{g})e^{-\frac{1}{2}t'\Psi_{n,\alpha}(\theta_{g})t}| \\
&\leq |\pi(\hat{\theta}_{n}+\frac{t}{\sqrt{n}})||\exp(z_{n})-\exp(z)|+ \exp(z)|\pi(\hat{\theta}_{n}+\frac{t}{\sqrt{n}})-\pi(\theta_{g})| \rightarrow 0
\end{aligned}
$$ as $n \rightarrow \infty$  by continuity of $\pi$ at $\theta_{g}$,  where $z_{n}=\exp(Q_{n}(\hat{\theta}_{n}+\frac{t}{\sqrt{n}}))$ and $z=e^{-\frac{1}{2}t'\Psi_{n,\alpha}(\theta_{g})t}$.
Now, for $t \in S_{2}$, using assumption (A5), choose $\delta_{0}$ so small that $|R_{n}(t)|<\dfrac{1}{4}t'\hat{\Psi}_{n}(\hat{\theta}_{n})t$, for all sufficiently large $n$. So for $t \in S_{2}$, we get 
$$
Q_{n}(\hat{\theta}_{n}+\frac{t}{\sqrt{n}})-Q_{n}(\hat{\theta}_{n})< -\dfrac{1}{4}t'\hat{\Psi}_{n}(\hat{\theta}_{n})t
$$
i.e.,
$$
\exp\left[Q_{n}(\hat{\theta}_{n}+\frac{t}{\sqrt{n}})-Q_{n}(\hat{\theta}_{n})\right]< e^{-\frac{1}{4}t'\hat{\Psi}_{n}(\hat{\theta}_{n})t}<e^{-\frac{1}{8}t'\Psi_{n,\alpha}(\theta_{g})t}.
$$
Hence, for $t \in S_{2}$,
$$|g_{n}(t)| \leq 2\pi(\theta_{g})e^{-\frac{1}{8}t'\Psi_{n,\alpha}(\theta_{g})t}+\pi(\theta_{g})e^{-\frac{1}{2}t'\Psi_{n,\alpha}(\theta_{g})t}
$$
 which is integrable. Hence by Dominated Convergence Theorem, 
$\int_{S_{2}}|g_{n}(t)|dt \rightarrow 0$ as $n \rightarrow \infty$.

Now we consider the integral over $S_{1}$. For $t \in S_{1}$, we get 
$$
\begin{aligned}
\dfrac{1}{n}\left[Q_{n}\left(\hat{\theta}_{n}+\frac{t}{\sqrt{n}} \right)-Q_{n}(\hat{\theta}_{n})\right]
 &=  \dfrac{1}{n}\left[Q_{n}\left(\hat{\theta}_{n}+\frac{t}{\sqrt{n}}\right)-Q_{n}(\theta_{g})\right]+\dfrac{1}{n}[Q_{n}(\theta_{g})-Q_{n}(\hat{\theta}_{n})]\\
&\leq  \sup\limits_{||\theta-\theta_{g}||>\frac{\delta_{0}}{2}}\frac{1}{n}(Q_{n}(\theta)-Q_{n}(\theta_{g}))\\
&+\dfrac{1}{2n}(\hat{\theta}_{n}-\theta_{g})'[\bigtriangledown^{2} Q_{n}(\hat{\theta}_{n})](\hat{\theta}_{n}-\theta_{g})\\
&+\dfrac{1}{6n}\sum_{j,k,l}(\hat{\theta}_{nj}-\theta_{gj})(\hat{\theta}_{nk}-\theta_{gk})(\hat{\theta}_{nl}-\theta_{gl})\bigtriangledown_{jkl}Q_{n}(\theta_{n}^{*})
\end{aligned}
$$
where $\theta_{n}^{*}$ lies between $\hat{\theta}_{n}$ and $\theta_{g}$. The first term in the last inequality comes from the
fact that $\hat{\theta}_{n}$ is consistent for $\theta_{g}$ and $\dfrac{||t||}{\sqrt{n}}>\delta_{0}$ as $t \in S_{1}$. Now using Assumption (E3), Assumption (A5) of Ghosh and Basu(2013) \cite{Ghosh/Basu:2013} and  $-\dfrac{1}{n}\bigtriangledown^{2} Q_{n}(\hat{\theta}_{n})=\hat{\Psi}_{n}(\hat{\theta}_{n})=\Psi_{n,\alpha}(\theta_{g})+o_{p}(1)$, it is clear that the second and third term of the above inequality goes to 0 with probability tending to 1 for sufficiently large $n$. By Assumption (E3), the first term is less than $-\epsilon$ with probability one for sufficiently large $n$ for some $\epsilon >0$, and hence with probability one, $\dfrac{1}{n}[Q_{n}\left(\hat{\theta}_{n}+\frac{t}{\sqrt{n}} \right)-Q_{n}(\hat{\theta}_{n})] <-\dfrac{\epsilon}{2}$, for all sufficiently large $n$. Therefore, we get\\

$\int_{S_{1}}|g_{n}(t)|dt \leq \int_{S_{1}}\pi(\hat{\theta}_{n}+\frac{t}{\sqrt{n}})e^{-\frac{n\epsilon}{2}}+\int_{S_{1}}\pi(\theta_{g})e^{-\frac{1}{2}t'\Psi_{n,\alpha}(\theta_{g})t}\\\\ \hspace*{0.92 in}\leq e^{-\frac{n\epsilon}{2}}\sqrt{n}\int \pi(\theta)d\theta+\pi(\theta_{g})\left(\dfrac{2\pi}{|\Psi_{n,\alpha}(\theta_{g})|}\right)^{\frac{p}{2}}\int\limits_{S_{1}}\phi_{n}(t)dt$\\\\
where $\phi_{n}(t)$ is the pdf of $N_{p}(\textbf{0},\Psi^{-1}_{n}(\theta_{g}))$. Now, the first integral goes to $0$, as $n \rightarrow \infty$ if $\int\pi(\theta)d\theta < \infty$, and using assumption (E3), the second integral also goes to 0, since this is a Normal Tail Probability, whose variance is non-zero and finite for sufficiently large $n$. Hence we have shown that, with probability tending to one, $\int\limits_{S_{1}}|g_{n}(t)| \rightarrow 0$, as $n \rightarrow \infty$, and that completes the proof of the first part.

The proof of the second part follows from the fact that  $\hat{\Psi}_{n}(\hat{\theta}_{n})=\Psi_{n,\alpha}(\theta_{g})+o_{p}(1)$.

\section{Verification of (E1)-(E2) for Fixed-design Logistic Regression}
\label{APP:logistic_verification}

Let $\beta^g$ be the best fitting parameter. So the true data generating density for $i^{th}$ observation is $g_i=f_i(x, \beta^g)$. 
Recall the notations of section 2. Note that $V_i(x,\beta)=\dfrac{1+e^{((1+\alpha)z_i^{T}\beta)}}{(1+e^{z_i^{T}\beta})^{1+\alpha}}-(1+\frac{1}{\alpha})\dfrac{e^{x\alpha z_i^{T}\beta}}{(1+e^{z_i^{T}\beta})^{\alpha}}$. 
Let $h_x(t)=\dfrac{e^{tx}}{1+e^t}$, $x=0,1$. Then $ h_x'(t)=\dfrac{xe^{tx}}{1+e^t}-\dfrac{e^{t(x+1)}}{(1+e^t)^2} $ ,  $h_x''(t)=x(\dfrac{xe^{tx}}{1+e^t}-\dfrac{e^{t(x+1)}}{(1+e^t)^2})-(\dfrac{(x+1)e^{t(x+1)}}{(1+e^t)^2}- 2\dfrac{e^{t(x+2)}}{(1+e^t)^3} )$, $h_x'''(t)=x^2(\dfrac{xe^{tx}}{1+e^t}-\dfrac{e^{t(x+1)}}{(1+e^t)^2})-(2x+1)\Big[(1+x)\dfrac{e^{t(x+1)}}{(1+e^t)^2}+2\dfrac{e^{t(x+2)}}{(1+e^t)^3} \Big] + 2 \Big[(x+2)\dfrac{e^{t(x+2)}}{(1+e^t)^3} +3 \dfrac{e^{t(x+3)}}{(1+e^t)^4} \Big]$.\\\\
So, $h_x'(t),h_x''(t),h_x'''(t)$ all are bounded for $x=0,1$.
 \\
Note that $$\nabla_jV_i(x,\beta)=z_{ij}(1+\alpha)\Big[\sum_{u=0}^{1}(h_u(z_i'\beta))^{\alpha}h_u'(z_i'\beta)-(h_x(z_i'\beta))^{\alpha-1}h_x'(z_i'\beta)\Big],$$
$$
\begin{aligned}
\nabla_{jk}V_i(x,\beta)=z_{ij}z_{ik}(1+\alpha)\Big[& \alpha \sum_{u=0}^{1}(h_u(z_i'\beta))^{\alpha-1}(h_u'(z_i'\beta))^2+\sum_{u=0}^{1}(h_u(z_i'\beta))^{\alpha}h_u''(z_i'\beta)\\ & -(\alpha-1)(h_y(z_i'\beta))^{\alpha-2}(h_x'(z_i'\beta))^2-(h_x(z_i'\beta))^{\alpha-1}h_x''(z_i'\beta)\Big],
\end{aligned}$$
$$
\begin{aligned}
\nabla_{jkl}V_i(x,\beta)=z_{ij}z_{ik}z_{il}(1+\alpha)\Big[& \alpha(\alpha-1) \sum_{u=0}^{1}(h_u(z_i'\beta))^{\alpha-2}(h_u'(z_i'\beta))^3+ 3\alpha \sum_{u=0}^{1}(h_u(z_i'\beta))^{\alpha-1}h_u'(z_i'\beta)h_u''(z_i'\beta)\\ & + \sum_{u=0}^{1}(h_u(z_i'\beta))^{\alpha}h_u'''(z_i'\beta) -(\alpha-1)(\alpha-2)(h_y(z_i'\beta))^{\alpha-3}(h_y'(z_i'\beta))^3 \\ & -3(\alpha-1)(h_x(z_i'\beta))^{\alpha-2}h_x'(z_i'\beta)h_x''(z_i'\beta)-(h_u(z_i'\beta))^{\alpha-1}h_x'''(z_i'\beta)\Big].
\end{aligned}$$

We first prove the condition (E1). So we have to prove the conditions (A1)-(A6) of Ghosh and Basu (2013) \cite{Ghosh/Basu:2013}. (A1)-(A3) are trivial, since all the model densities are thrice differentiable, and we can take the differential under the integral sign. (A4) is obvious from condition L2, since the condition (A4) is (R3) itself. 
Now we have to show (A5), i.e there exists a function $M_{jkl}^{(i)}(x)$ such that 
$$|\nabla_{jkl}V_i(x,\beta)| \leq M_{jkl}^{(i)}(x) \hspace*{0.1 in} \forall \beta, \forall i $$
where 
$$\dfrac{1}{n}\sum_{i=1}^{n}E_{g_i}\Big[M_{jkl}^{(i)}(X)\Big]=O(1)\hspace*{0.2 in} \forall j,k,l.$$
Now R1 and the above expression for $\nabla_{jkl}V_i(x,\beta)$ implies $|\nabla_{jkl}V_i(x,\beta)| \leq C$ for some constant C, which obviously imply (A5). \\
Now we focus on (A6). We have to show 
$$\lim\limits_{N \rightarrow \infty} \sup\limits_{n >1}\Big\{ \dfrac{1}{n} \sum\limits_{i=1}^{n} E_{g_{i}}[|\nabla_{j}V_{i}(X,\theta)|I(|\nabla_{j}V_{i}(X,\theta)|>N)]\Big\}=0$$ and 
$$
\begin{aligned}
\lim\limits_{N \rightarrow \infty} \sup\limits_{n >1}\Big\{ \dfrac{1}{n} \sum\limits_{i=1}^{n} & E_{g_{i}}[|\nabla_{jk}V_{i}(X,\theta)-E_{g_{i}}(\nabla_{jk}V_{i}(X,\theta))|\\ & I(|\nabla_{jk}V_{i}(X,\theta)-E_{g_{i}}(\nabla_{jk}V_{i}(X,\theta))|>N)]\Big\}=0.
\end{aligned}$$
Since R1 says that $ \displaystyle{\max _{1 \leq i \leq n }} |z_{ij}||z_{ik}|=O(1) $, hence $\nabla_{jk}V_{i}(x,\theta)$ is bounded for $x=0,1$. Hence for large $N$, the indicator itself becomes 0, regardless of $i,j,k$. This finally proves (E1).\\
Now come to E2. Throughout the rest of the proof, we write $Q_n^{(\alpha)}(\beta)$ as $Q_n^(\beta)$ for notational simplicity. Note that 
$$Q_n(\beta)=\dfrac{1}{\alpha}\sum_{i=1}^{n}e^{\alpha(x_iz_i'\beta-\log(1+e^{z_i'\beta}))}-\dfrac{1}{1+\alpha}\sum_{i=1}^{n}\Big[e^{(1+\alpha)(-\log(1+e^{z_i'\beta}))}+e^{(1+\alpha)(z_i'\beta-\log(1+e^{z_i'\beta}))}\Big].$$
So, using MVT on each of the term in the summand for $Q_n(\beta)$, and write them as a function of $z_i'\beta$, we get
$$\dfrac{1}{n}(Q_n(\beta)-Q_n(\beta^g))=\dfrac{1}{n}(\beta-\beta^g)'\sum_{i=1}^{n}\dfrac{e^{\alpha t_i^{*}x_i}}{(1+e^{t_i^{*}})^{1+\alpha}}(x_i-\dfrac{e^{t_i^{*}}}{1+e^{t_i^{*}}})z_i-\dfrac{1}{n}(\beta-\beta^g)'\sum_{i=1}^{n}\dfrac{e^{(1+\alpha)t_i^{*}}-e^{t_i^{*}}}{(1+e^{t_i^{*}})^{2+\alpha}}z_i$$
where $t_i^{*}$ is between $z_i'\beta$ and $z_i'\beta^g$.
Since $x_i$ can be only 0 and 1, so write $$e^{\alpha t_i^{*}x_i}=x_ie^{\alpha t_i^{*}}+(1-x_i).$$
Hence, after some calculation.
$$\dfrac{1}{n}(Q_n(\beta)-Q_n(\beta^g))=\dfrac{1}{n}(\beta-\beta^g)' \sum_{i=1}^{n}\dfrac{x_i(e^{\alpha t_i^{*} }+e^{t_i^{*}})}{(1+e^{t_i^{*}})^{1+\alpha}}z_i-\dfrac{1}{n}(\beta-\beta^g)' \sum_{i=1}^{n}\dfrac{(e^{(1+\alpha) t_i^{*} }+e^{2t_i^{*}})}{(1+e^{t_i^{*}})^{2+\alpha}}z_i.
$$
Now $E_{\beta^g}[x_i]=\dfrac{e^{z_i'\beta^g}}{1+e^{z_i'\beta^g}}$. Since all terms in the coefficient of $x_i$ in the above term is bounded, so using Kolmogorov SLLN, we get 
$$v_n=\dfrac{1}{n}\sum_{i=1}^{n}\dfrac{x_i(e^{\alpha t_i^{*} }+e^{t_i^{*}})}{(1+e^{t_i^{*}})^{1+\alpha}}z_i-\dfrac{1}{n}\sum_{i=1}^{n}\dfrac{e^{z_i'\beta^g}}{1+e^{z_i'\beta^g}}\dfrac{(e^{\alpha t_i^{*} }+e^{t_i^{*}})}{(1+e^{t_i^{*}})^{1+\alpha}}z_i \xrightarrow[]{a.s} 0.$$
Hence,
$$
\begin{aligned}
& \dfrac{1}{n}(Q_n(\beta)-Q_n(\beta^g))=-\dfrac{1}{n}(\beta-\beta^g)' \sum_{i=1}^{n}(\dfrac{e^{t_i^{*}}}{1+e^{t_i^{*}}}-\dfrac{e^{z_i'\beta^g}}{1+e^{z_i'\beta^g}})\dfrac{(e^{\alpha t_i^{*} }+e^{t_i^{*}})}{(1+e^{t_i^{*}})^{1+\alpha}}z_i+(\beta-\beta^g)'v_n\\ & 
=-\dfrac{1}{n}(\beta-\beta^g)' \sum_{i=1}^{n}\dfrac{e^{t_i^{**}}}{(1+e^{t_i^{**}})^2}\dfrac{(e^{\alpha t_i^{*} }+e^{t_i^{*}})}{(1+e^{t_i^{*}})^{1+\alpha}}z_i(t_i^{*}-z_i'\beta^g)+(\beta-\beta^g)'v_n
\end{aligned}
$$
where $t_i^{**}$ is between $ t_i^{*} $ and $ z_i'\beta^g $. Note that $(z_i'\beta-z_i'\beta^g)(t_i^{*}-z_i'\beta^g) \geq (t_i^{*}-z_i'\beta^g)^2.$\\
Now, we can get a $\beta^{q} \neq \beta^g$ such that $=(t_i^{*}-z_i'\beta^g)^2 \geq (\beta^{q}-\beta^g)'z_iz_i'(\beta^{q}-\beta^g)$ for all $\beta$ outside the $\delta$ neighbourhood of $\beta^g$ . So 
$$\dfrac{1}{n}(\beta-\beta^g)' \sum_{i=1}^{n}\dfrac{e^{t_i^{**}}}{(1+e^{t_i^{**}})^2}\dfrac{(e^{\alpha t_i^{*} }+e^{t_i^{*}})}{(1+e^{t_i^{*}})^{1+\alpha}}z_i(t_i^{*}-z_i'\beta^g) \geq \dfrac{1}{n}(\beta^{q}-\beta^g)' \sum_{i=1}^{n}\dfrac{e^{t_i^{**}}}{(1+e^{t_i^{**}})^2}\dfrac{(e^{\alpha t_i^{*} }+e^{t_i^{*}})}{(1+e^{t_i^{*}})^{1+\alpha}}z_iz_i'(\beta^{q}-\beta^g) .$$
Note that, by supporting hyperplane theorem, $\exists$ $ \beta_0$ such that $z_i'\beta_0 > 0 $ for all $i=1,2,...,n$. Also note that the two functions $h_1(t)=\dfrac{e^{t}}{(1+e^{t})^2}$ and $h_2(t)=\dfrac{e^{\alpha t }+e^{t}}{(1+e^{t})^{1+\alpha}}$ are even functions, i.e $h_1(t)=h_1(-t)$ and $h_2(t)=h_2(-t)$ .  So, w.l.o.g, we assume $t_i^{*},t_i^{**}>0$. Since $h_1(t)$ is decreasing for $t>0$, and $h_2(t) $ is decreasing for $t>\frac{1}{\alpha}$, hence we can find a positive constant $C$ such that $h_1(t_i^{**})\geq h_1(z_i'\beta^g+Cz_i'\beta_0)  $, and $h_2(t_i^{*})\geq h_2(z_i'\beta^g+Cz_i'\beta_0)$. Call $\beta^{*}=\beta^g+C\beta_0$. So finally we get, 
$$
\begin{aligned}
\dfrac{1}{n}(Q_n(\beta)-Q_n(\beta^g)) & \leq -\dfrac{1}{n}(\beta^q-\beta^g)' \sum_{i=1}^{n} \dfrac{e^{z_i^{T}\beta^{*}}(e^{\alpha z_i^{T}\beta^{*}}+e^{z_i^{T}\beta^{*}})}{(1+e^{z_i^{T}\beta^{*}})^{3+\alpha}}z_iz_i^{T} (\beta-\beta^g) + (\beta^q-\beta^g)'v_n\\
& =-(\beta^q-\beta^g)'\Psi_n(\beta^{*})(\beta^q-\beta^g)+ (\beta-\beta^g)'v_n.
\end{aligned}$$
Now using R2 and same argument as we have proved the conditions in linear regression case, $(\beta^q-\beta^g)'\Psi_n(\beta^{*})(\beta^q-\beta^g) > \eta \delta^2 $ for some $\eta>0$ on the region $||\beta-\beta^g||>\delta$. Take $\epsilon=\dfrac{\eta \delta^2}{2}$. Then E2 holds for this specific $\epsilon$. E3 directly follows from R3.


\section{Proof of Theorem 5.2}\label{APP:BP}

Since we fix a scale parameter, let us denote the modelling density by $f_{i,\mu}$. Note that the ERPE functional is given by $$ T_{n}^{(\alpha)e}(\mathcal{G})=\dfrac{\int \mu \exp[Q^{(\alpha)}(\mu; \mathcal{G},\mathcal{F}_{n,\mu})]\pi(\mu)d\mu}{\int \exp[Q^{(\alpha)}(\mu; \mathcal{G},\mathcal{F}_{n,\mu})]\pi(\mu)d\mu}.$$
Consider the contamination distributions  $\{K_{i,m}\}$. Denote $H_{1,m}=(1-\epsilon)G_1+\epsilon K_{1,m},..., H_{n,m}=(1-\epsilon)G_n+\epsilon K_{n,m}$, and subsequently $\mathcal{H}_m=(H_{1,m},...,H_{n,m})$.
By Theorem \ref{THM:Expansion_Func}
\begin{equation}
T_{n}^{(\alpha)e}(\mathcal{G})\sim T_{n,\alpha}(\mathcal{G})\Big[1+O(n^{-1})\Big],
\end{equation} 
and
 $$T_{n}^{(\alpha)e}(\mathcal{H}_m)\sim T_{n,\alpha}(\mathcal{H}_m)\Big[1+O(n^{-1})\Big],$$
 where $T_{n,\alpha}$ is the minimum DPD functional. Hence we clearly see that if $\lvert T_{n,\alpha}(\mathcal{G})- T_{n,\alpha}(\mathcal{H}_m) \rvert$ doesn't go to $\infty$ for some $\epsilon>0$ as $m\rightarrow \infty$, $\lvert T_{n}^{(\alpha)e}(\mathcal{G})- T_{n}^{(\alpha)e}(\mathcal{H}_m) \rvert$ also doesn't go to $\infty$ asymptotically. Hence the asymptotic breakdown point for  $ T_{n}^{(\alpha)e} $ is at least $\frac{1}{2}$.



\end{document}